\documentclass[lettersize,journal]{IEEEtran}

\usepackage{cite}
\usepackage{amsmath,amssymb,amsfonts}
\usepackage{algorithmic}
\usepackage{graphicx}
\usepackage{textcomp}
\usepackage{mathtools}
\usepackage{subfig}
\usepackage{float}
\usepackage{graphicx}
\usepackage{amssymb}
\usepackage{epsfig} % for postscript graphics files
\usepackage{algorithm}
\usepackage{stfloats}
\newtheorem{theorem}{Theorem}[section]

\newtheorem{proposition}[theorem]{Proposition}

\newtheorem{definition}[theorem]{Definition}
\newtheorem{remark}[theorem]{Remark}
\newtheorem{assumption}[theorem]{Assumption}
\DeclareMathOperator{\tr}{tr}
\newcommand\blfootnote[1]{%
  \begingroup
  \renewcommand\thefootnote{}\footnote{#1}%
  \addtocounter{footnote}{-1}%
  \endgroup
}

\begin{document}

\textbf{This work has been submitted to the IEEE for possible publication. Copyright may be transferred without notice, after which this version may no longer be accessible.}
\title{Reinforcement Learning for Inverse Non-Cooperative Linear-Quadratic Output-feedback Differential Games}

\author{Emin Martirosyan, Ming Cao\\
Engineering and Technology Institute Groningen\\ University of Groningen, Nijenborgh 4, Groningen, 9712CP, Netherlands}

% The paper headers
\markboth{Journal of \LaTeX\ Class Files,~Vol.~14, No.~8, August~2021}%
{Shell \MakeLowercase{\textit{et al.}}: A Sample Article Using IEEEtran.cls for IEEE Journals}

\maketitle
\blfootnote{\textit{\textbf{This is a preprint of work that has been submitted to the IEEE for possible publication. Copyright may be transferred without notice, after which this version may no longer be accessible..}}}

%%%%%%%%%%%%%%%%%%%%%%%%%%%%%%%%%%%%%%%%%%%%%%%%%%%%%%%%%%%%%%%%%%%%%%%%%%%%%%%%
\begin{abstract}
In this paper, we address the inverse problem for linear-quadratic differential non-cooperative games with output-feedback. Given players' stabilizing feedback laws, the goal is to find cost function parameters that lead to a game for which the observed game dynamics are at a Nash equilibrium. Using the given feedback laws, we introduce a model-based algorithm that generates cost function parameters solving the above inverse problem. We introduce a correction procedure that at each iteration of the algorithm guarantees the existence of the feedback laws, which addresses a key challenge of output-feedback control designs. As an intermediate stage of the algorithm, we have developed a procedure for the initial stabilization of the multiple-input system with output-feedback information structure. We prove convergence and stability of the algorithm, and show the way to generate new games with necessary properties without requiring to run the complete algorithm repeatedly. Then the algorithm is extended to a model-free version that uses data samples generated by unknown dynamics and has the same converging and stabilizing properties as the model-based version. Finally, we show how the inverse problem can be solved in a distributed manner and provide possible extensions. Simulation results validate the effectiveness of the proposed algorithms.
\end{abstract}


\section{Introduction}
\label{sec:introduction}

\IEEEPARstart{G}ame theory examines multi-person decision-making mechanisms. A game becomes dynamic when decision order matters and players influence a common dynamical system \cite{molloy2022inverse}. If each player pursues their own interests, often conflicting with others', the game is non-cooperative \cite{bacsar1998dynamic}. The fusion of dynamic and non-cooperative elements has drawn attention from the control community, offering fresh perspectives on control law performance evaluation. Initially proposed by Isaacs in the 1960s for two-player non-cooperative pursuit-evasion games \cite{isaacs1965differential}, dynamic non-cooperative games were further developed in the seminal work by Starr \cite{starr1969nonzero}, introducing $N$-player variants. Subsequent research applied dynamic games to vehicle collision avoidance \cite{molloy2020optimal}, decentralized control of electric vehicles \cite{ma2011decentralized}, vehicle formation control \cite{gu2007differential}, advanced driver assistance systems \cite{na2014game}, and human-machine interaction or shared control systems modeling \cite{inga2021online}. Existing literature predominantly focuses on determining game outcomes and player behavior based on objective functions, with only recent attention directed toward inverse problems. Inverse dynamic games aim to reverse engineer player objectives from observed gameplay behavior.

The evolution of non-cooperative dynamic game theory has largely paralleled optimal control, yet surprisingly, until recently, little attention has been devoted to its inverse problem \cite{molloy2022inverse}. Inverse dynamic games, a generalization of inverse optimal control (IOC), broaden the scope from a single dynamical system and objective function to potentially multiple systems and objectives \cite{tsai2016inverse}. IOC, dating back to 1966 \cite{anderson1966inverse}, focuses on inferring system objectives and constraints from observed behavior. Inverse reinforcement learning (IRL) \cite{ng2000algorithms} tackles a similar task within the framework of Markov decision processes. Efforts on inverse dynamic games include solving finite-time LQ in \cite{kopf2017inverse}, utilizing an extended maximum entropy IRL \cite{ziebart2008maximum}. Other works like \cite{molloy2019inverse,rothfuss2017inverse} establish connections between inverse optimal control and inverse dynamic games, providing solutions to inverse open-loop differential games. Further research has addressed imitation problems/apprentice games \cite{lian2022data, lian2021inverse, lian2022inverse, lian2023off}, akin to inverse dynamic game problems. Dynamic games' close relation to optimal control warrants exploration of literature dedicated to tracking and inverse optimal control, typically focused on single-input systems or single optimization objectives \cite{xue2021inversec,xue2021inverse}. Most works assume player control inputs formed from system state, except for \cite{xue2023inverse} and \cite{lian2023inverse}, which consider zero-sum differential games with output-feedback and discrete-time linear systems with output-feedback information structures, respectively. The algorithm introduced typically yields sub-optimal behavior for the designed game.

In this work, we consider games with linear continuous-time system dynamics and quadratic cost functions, where players, instead of using the state dynamics to form their control inputs, make use of their own heterogeneous outputs. To deal with the output-feedback information structure, we used results developed in \cite{ilka2022novel,engwerda2008result} for the linear quadratic regulator (LQR) problem. Firstly, a model-based algorithm that finds the cost function for all the players is developed. The desired behavior is Nash optimal for the designed game. Then, the algorithm is extended to the model-free setting. Using ideas presented in \cite{modares__2015,jiang2012computational}, and the so-called integral reinforcement learning \cite{vrabie2009adaptive}, we assume that instead of the system matrices, some trajectories generated by the unknown linear system can be observed. We also characterize all possible solutions for the inverse differential game with the described players' information structure. Finally, we demonstrate how the inverse differential game problem can be solved in a distributed manner.

The paper is structured as follows. Section \ref{sec:ProblemFormulation} shows preliminary results on LQ non-cooperative $N$-player differential games with output-feedback, and formulates the problem addressed in this paper. In section \ref{sec:MB}, we describe each step of the model-based algorithm, present its analytical properties, and characterize possible solutions of the inverse problem. Section \ref{sec:ModelFree} extends the model-based algorithm to a model-free version that allows to solve the problem without using the system matrices, and introduces the distributed procedure that allows to solve the inverse problem. In section \ref{sec:Siumation}, we provide simulation results that validate the effectiveness of the proposed algorithms. Finally, section \ref{sec:Conclusion} concludes the paper with suggestions on the possible future research. 

\textit{Notations}: For a matrix $P\in\mathbb{R}^{m\times n}$, $P^k$ and $P^{(k)}$ denote $P$ to the power of $k$, and the matrix $P$ at the $k$-th iteration, respectively. In addition, $P>0$ and $P\geq0$ denote positive definiteness and semi-definitness of the matrix $P$, respectively. The notation $\tr P$ denotes the trace of the matrix $P$. $I_k$ and $\mathbf{0}_k$ are the $k\times k$ identity matrix and zero matrices, respectively.

\section{Problem Formulation}
\label{sec:ProblemFormulation}
In this section, LQ output-feedback differential games are introduced. We define the stationary linear feedback Nash equilibrium (further referred to as NE) and introduce the inverse problem for the given class of games. 

Consider a differential game with $N$ players, where we denote the set of players by $\mathcal{N} = \{1,\dots,N\}$, under the continuous time dynamics of the game 
\begin{align}\label{sysdyn}
    \begin{split}
        \dot{x}(t) &= A x(t) + \sum_{i=1}^N B_i u_i(t),\quad x(0) = x_0, \\
        y_i(t) &= C_i x(t),
    \end{split}
\end{align}
where $x\in\mathbb{R}^N$ is the state, $u_i\in\mathbb{R}^{m_i}$ and $y_i\in\mathbb{R}^{p_i}$ ($p_i < n$) are the control input and the output of player $i$, respectively; the plant matrix $A$, control input matrices $B_i$ and output matrices $C_i$ have appropriate dimensions.

We consider that the players select their control input to be 
\begin{equation}\label{continputs}
    u_i(t) = - K_i y_i(t),
\end{equation}
where $K_i$ is an $m_i\times p_i$ time-invariant feedback matrix of player $i$. Further, to ease notations, we use $x(t) = x$, $u_i(t)=x$ and $y_i(t) = y_i$; $F_i = K_i C_i$  for $i \in\mathcal{N}$.

\begin{assumption}\label{stabassumption}
    System \eqref{sysdyn} is output-feedback stabilizable, i.e, there exists a tuple $(K_1,\dots,K_N)$ such that 
    \begin{equation}
        A - \sum_{i=1}^N B_i K_i C_i 
    \end{equation}
    is stable. 
\end{assumption}
\begin{assumption}
    $C_i$ is full row rank for $i\in\mathcal{N}$.
\end{assumption}

We use $K_{-i} = (K_1,\dots,K_{i-1},K_{i+1},\dots,K_N)$ to denote a profile of feedback laws of all the players except for player $i$. Within the game each of the players aims to minimize its cost function $J_i(x_0,u_i,u_{-i})$ which has the quadratic form 
\begin{equation}\label{costfun}
  J_i(x_0,u_i,u_{-i}) = \int_t^\infty \big( x^\top Q_i x + \sum_{j=1}^N u_j^\top R_{ij} u_j \big)\,d\tau,
\end{equation}
where $Q_i$ and $R_{ij}$ are symmetric matrices, and $R_{ii} > 0$.
Players are considered to choose their feedback laws from the following set 
\begin{equation}
    \mathcal{K} = \{(K_1,\dots,K_N) | A + \sum_{i=1}^N B_i K_i C_i\quad\text{is stable}\}.
\end{equation}
This restriction spoils the rectangular structure of the strategy spaces, i.e. choices of feedback matrices cannot be made independently. However, such a restriction is motivated by the fact that closed-loop stability is usually a common objective \cite{engwerda2008result}. Assumption \ref{stabassumption} guarantees the non-emptiness of $\mathcal{K}$.

Further, we provide the definition of output-feedback NE introduced in \cite{engwerda2008result}. 
\begin{definition}
    An $N$-tuple $(K_1^*,\dots,K_N^*)\in\mathcal{K}$ is called an output-feedback NE if for every player $i\in\mathcal{N}$ the following inequality holds 
    \begin{equation}
        J_i(x_0, K_i^*, K_{-i}^*) \leq J_i(x_0, K_i, K_{-i}^*)
    \end{equation}
    for each $x_0$ and for each $K_i$ such that $(K_1^*,\dots,K_{i-1}^*,K_i,K_{i+1}^*,\dots,K_N^*)\in\mathcal{K}$.
\end{definition}

Next, let us introduce the set of algebraic Riccati equations (AREs)
\begin{align}\label{ARE}
    \begin{split}
        & A^\top X_i + X_i  A - \sum_{j\neq i}^N X_i B_j R_{jj}^{-1} B_j^\top X_j - \\ & \sum_{j\neq i}^N X_j B_j R_{jj}^{-1} B_j^\top X_i - X_i B_i R_{ii}^{-1} B_i^\top X_i +  \\ 
        & \sum_{j\neq i}^N X_j B_j R_{jj}^{-1} R_{ij} R_{jj}^{-1} B_j^\top X_j + Q_i = 0.
    \end{split}
\end{align}
for $i\in\mathcal{N}$. The stabilizing solution of the above system of AREs is a tuple $(X_1,\dots,X_N)$ such that $A - \sum_{u=1}^N B_i R_{ii}^{-1} B_i^\top X_i$ is stable. Next, we introduce the main known theoretical result presented in \cite{engwerda2008result} that relates the output-feedback NE in the described games and stabilizing solutions of \eqref{ARE}. Further, $X_i$ is referred to as a value matrix for $i\in\mathcal{N}$.

\begin{theorem}\label{mainth}
    Let $(X_1,\dots,X_N)$ be a stabilizing solution of \eqref{ARE} and $(K_1,\dots,K_N)$ be such that $R_{ii}^{-1} B_i^\top X_i = K_i C_i$ for $i\in\mathcal{N}$. Then $(K_1,\dots,K_N)$ is an output-feedback NE. Conversely, if $(K_1,\dots,K_N)$ is an output-feedback NE, there exists a stabilizing solution $(X_1,\dots,X_N)$ of \eqref{ARE} and $(K_1,\dots,K_N)$ such that $R_{ii}^{-1} B_i^\top X_i = K_i C_i$. 

    Furthermore, if the game has an output-feedback Nash equilibrium with this equilibrium, the corresponding cost for player $i$ is $J_i (x_0, K_i, K_{-i}) =  x_0^\top X_i x_0$.
\end{theorem}

\begin{remark}
    Note that it is not enough to find a stabilizing solution $(X_1,\dots,X_N)$ for \eqref{ARE}. The solution should be such that for the given $(C_1,\dots,C_N)$ there exists $(K_1,\dots,K_N)$ satisfying $K_i C_i = R_{ii}^{-1} B_i^\top X_i$ for $i\in\mathcal{N}$. 
\end{remark}

We introduce the following notations 
    \begin{align}
            &(B_1,\dots,B_N) = \mathbf{B},\, (C_1,\dots,C_N) = \mathbf{C},\,(Q_1,\dots,Q_N) = \mathbf{Q},\nonumber\\
            &(R_{11},\dots,R_{1N},\dots,R_{N1},\dots,R_{NN}) = \mathbf{R},
    \end{align}
A tuple $(A,\mathbf{B},\mathbf{C},\mathbf{Q},\mathbf{R})$ refers to a game with dynamics \eqref{sysdyn} and cost functions parameters \eqref{costfun}. Now, for the given setting, we are ready to introduce the inverse problem for differential games.\newline
\textbf{Problem 1}. Suppose there is a tuple $(K_{1,d},\dots,K_{N,d})$ of controllers that stabilize a system described by $(A,\mathbf{B},\mathbf{C})$. Given matrices $(A,\mathbf{B},\mathbf{C})$, we want to construct quadratic cost functions as in \eqref{costfun} such that $(K_{1,d},\dots,K_{N,d})$ is an NE. \\
\textbf{Problem 2}. Suppose there is a tuple $(K_{1,d},\dots,K_{N,d})$ of controllers that stabilize a system described by $(A,\mathbf{B},\mathbf{C})$. Given matrices $\mathbf{C}$ and the feasibility to apply control inputs to the linear system described by matrices $(A,\mathbf{B},\mathbf{C})$ for collecting the data on resulted trajectories, we want to construct quadratic cost functions as in \eqref{costfun} such that $(K_{1,d},\dots,K_{N,d})$ is an NE. 

\textbf{Problem 1} and \textbf{Problem 2} are solved iteratively via the so-called model-based and model-free algorithms, respectively.
\begin{remark}
    Note that Assumption \ref{stabassumption} is not necessary for the given system because $(K_{1,d},\dots,K_{N,d})$ is known to be a stabilizable tuple which implies that the system is stabilizable. 
\end{remark}
However, another assumption needs to be done. 
\begin{assumption}\label{CBassump}
    It is assumed that $C_i B_i \neq \mathbf{0}_{p_i\times m_i}$ where $\mathbf{0}_{p_i\times m_i}$ is the zero matrix of dimension $p_i\times m_i$ for $i\in\mathcal{N}$.
\end{assumption}

The necessity for this assumption as well as the way to relax it are discussed in Section \ref{sec:Evaluation of the Parameters}. 
\begin{remark}
    Note that in \textbf{Problem 2}, matrices $(A,\mathbf{B})$ are assumed to be unknown. For this problem, we also show a way to solve the problem in a distributed way which is described in Section \ref{sec:Distributed Procedure}. 
\end{remark}
\begin{remark}
    In fact, the problem can be reformulated if instead of a tuple of desired controllers $(K_{1,d},\dots,K_{N,d})$, we are given a set of trajectories $\{(y_1,u_1),\dots,(y_N,u_N)\}$ generated by the this tuple of controllers \cite{inga2019solution}. Then, one can estimate $(K_{1,d},\dots,K_{N,d})$ from the given trajectories from \eqref{continputs}, e.g. using the least-square method \cite{devore2008probability}.  
\end{remark}
\begin{remark}\label{multiplesolutions}
   It is known that one tuple of target feedback laws can be optimal for different sets of the cost parameters \cite{lancaster1995algebraic}. In the following sections, we show how after derivation of a set of parameters with desired properties via the proposed algorithms, one can generate a new set not requiring using the algorithms again.  
\end{remark}

\section{Model-based Algorithm}
\label{sec:MB}
This section is dedicated solely to the development of the model-based algorithm. The whole procedure can be described in the following way - after the initialization of the parameters, we update the value matrices in the direction of desired feedback laws and consequently update the cost function parameters. Note that during the iterative procedure, matrices $R_{ij}$ remain the same while only $Q_i$ matrices are updated.  

The algorithm is designed in such a way that at each iteration, the feedback laws tuple, updated in the direction of the target tuple $(K_{1,d},\dots,K_{N,d})$, stabilizes the system and is Nash optimal for the cost function with parameters generated at the same iteration. This is a useful property for the tracking/imitation problems when one wants to keep the system stable during the tracking/imitation procedure \cite{lian2021inverse, lian2022data, lian2022inverse, lian2023off, xue2023inverse, xue2021inversec, xue2021inverse}. We show the procedure to find such $X_i$ for $i\in\mathcal{N}$ that guarantee the existence of $K_i$ for $i\in\mathcal{N}$ and stabilize the system. 

\subsection{Existence of the output-feedback}
Let us start from discussing the existence of $K_i$ for $X_i$ computed from \eqref{ARE} for $i\in\mathcal{N}$. Suppose we computed $X_i$ from \eqref{ARE}. According to \cite{engwerda2008result}, the output-feedback law $K_i$ exists if
\begin{equation}\label{existcond}
    R_{ii}^{-1} B_i^\top X_i = R_{ii}^{-1} B_i^\top X_i C_i^\top (C_i C_i^\top)^{-1} C_i = K_i C_i,
\end{equation}
where $R_{ii} > 0$. Thus, it is not enough only to compute the solution of \eqref{ARE}, and the following
\begin{equation}
    B_i^\top X_i (I_n -  C_i^\top (C_i C_i^\top)^{-1} C_i) = 0
\end{equation}
needs to be satisfied. From this point we define $C_i^+ =  C_i^\top (C_i C_i^\top)^{-1}$. 

Suppose one gets $X_i$ such that for some $K_i$ the following holds 
\begin{equation}
    K_i C_i = R_{ii}^{-1} B_i^\top X_i C_i^+ C_i\neq R_{ii}^{-1} B_i^\top X_i.
\end{equation}
Using the gradient descent optimization \cite{Bertsekas/99}, we can find $\tilde{X}_i$ so that $R_{ii}^{-1} B_i^\top \tilde{X}_i C_i^+ C_i = R_{ii}^{-1} B_i^\top \tilde{X}_i = K_i C_i$. Define 
\begin{equation}
    d_i^{(k)} \coloneq B_i^\top X_i^{(k)} - R_{ii} K_i C_i,
\end{equation}
where $k = 0,1,\dots$ is the iteration counter and $X_i^{(0)}$ is such that $R_{ii}^{-1} B_i^\top X_i^{(0)} C_i^+ = K_i$. We construct $d_i^{(k)\top} d_i^{(k)}$ which is convex in $X_i$ and has its minimum when $d_i^{(k)} = 0$. One can use the gradient descent update 
\begin{align}\label{existproc}
    \begin{split}
        X_i^{(k+1)}  & = X_i^{(k)} - \alpha_i \frac{\partial \tr (d_i^{(k)\top} d_i^{(k)})}{\partial X_i}\\
        & = X_i^{(k)} - \alpha_i (B_i d_i^{(k)} + d_i^{(k)\top} B_i^\top), \\
        %& d_i^{(k+1)} = B_i^\top X_i^{(k+1)} - R_{ii} K_i C_i,
    \end{split}
\end{align}
where $\alpha_i > 0$ is the step size.%, and the partial derivative is computed as 
% \begin{equation}
%      \frac{\partial (\tr (d_i^{(k)\top} d_i^{(k)})}{\partial X_i} = B_i d_i^{(k)} + d_i^{(k)\top} B_i^\top.
% \end{equation}

Following this procedure, one can find $\tilde{X}_i$ satisfying \eqref{existcond}, i.e., $R_{ii}^{-1} B_i^\top \tilde{X}_i = K_i C_i$. Further, we refer to the described gradient descent procedure as $\mathcal{G}(X_i, K_i,\alpha_i)$ and write 
\begin{equation}\label{corrfun}
\mathcal{G}(X_i, K_i,\alpha_i) = \tilde{X}_i    
\end{equation}
which means that $\tilde{X}_i$ satisfying $R_{ii}^{-1} B_i^\top \tilde{X}_i = K_i C_i$ was found via \eqref{existproc} with $\alpha_i$ step size and the initial value $X_i^{(0)} = X_i$ such that $R_{ii}^{-1} B_i^\top X_i^{(0)} C_i^+ = K_i$. 

\subsection{Initialization}
\label{sec:Initialization}

In this section, we initialize the parameters for \textbf{Problem 1}. To be more precise, for some fixed $R_{ii} > 0$ with $i,j \in\mathcal{N}$, we compute a tuple $(K_1^{(0)},\dots, K_N^{(0)})$ that is stabilizing, i.e., $A - \sum_{i=1}^N B_i K_i^{(0)} C_i$ is stable, where $K_i^{(0)} C_i = R_{ii}^{-1} B_i^\top X_i^{(0)}$. 

Finding a stabilizing output-feedback solution for linear time-invariant systems is mostly studied within the linear quadratic regulator (LQR) framework and it is a challenging problem that has been studied a lot \cite{syrmos1997static}. Our problem is simpler - we do not need to find the initial NE (optimal solution in the case of LQR) for the given set of the parameters, but need to find a tuple $(X_{1}^{(0)},\dots,X_{N}^{(0)})$ that forms stabilizing $(K_1^{(0)},\dots,K_N^{(0)})$. %and then we can adjust $Q_i^{(0)}$ for each player $i\in\mathcal{N}$.
Towards this goal, we extend the algorithm designed for LQR with the output-feedback presented in \cite{ilka2022novel} (Algorithm 2: Modified Newton's Method for Static output-feedback Controller Design). 

The first step of the extended algorithm is the computation of some initialized tuple $(P_1,\dots,P_N)$ such that $A - \sum_{i=1}^N B_i R_{ii}^{-1} B_i^\top P_i$ is stable. We find a stabilizing tuple $(P_1,\dots,P_N)$ by solving the set AREs associated with state-feedback differential games. In \cite{engwerda2007algorithms}, set AREs are solved iteratively until the convergence of $P_i$ for $i\in\mathcal{N}$ to derive NE for the state-feedback differential games. However, in our case, we do not need the NE solution but just a stabilizing one. Hence, the following AREs
\begin{align}\label{lyapiter}
    \begin{split}
        & (A - \sum_{j=i}^N B_j F_j - \sum_{j=1}^{i-1} B_j R_{jj}^{-1} B_j^\top P_j)^\top P_i + \\
        & P_i (A - \sum_{j=i}^N B_j F_j - \sum_{j=1}^{i-1} B_j R_{jj}^{-1} B_j^\top P_j) =\\
        & - (C_i^\top C_i - F_i^\top R_{ii} F_i)
    \end{split}
\end{align}
are solved asynchronously (first for $i=1$, then $i=2$, etc.) only once with respect to $P_i$ for each $i\in\mathcal{N}$, given a stable $A - \sum_{j=1}^N B_j F_j$. The result is a stabilizing tuple $(P_1,\dots,P_N)$, i.e., the matrix
\begin{equation}
     A - \sum_{j=1}^N B_j R_{jj}^{-1} B_j^\top P_j
\end{equation}
is stable. To guarantee that $A - \sum_{j=1}^N B_j F_j$ is stable, $F_j$ can be set to be $K_{j,d} C_j$, i.e., $F_j = F_{j,d} = K_{j,d} C_j$ for $j\in\mathcal{N}$ that are known to satisfy
\begin{equation}
    A - \sum_{j=1}^N B_j K_{j,d} C_j\quad\text{is stable}.
\end{equation}
In fact, $F_j$ for $j\in\mathcal{N}$ can be set to be any stabilizing tuple $(F_1,\dots,F_N)$ such that $A - \sum_{j=1}^N B_j F_j$ is stable. However, because from the problem formulation, we already have given the NE tuple $(K_{1,d},\dots,K_{N,d})$, we set $F_j = F_{j,d}$.

Using the derived stabilizing tuple of $(P_1,\dots,P_N)$, we can extend the modified Newton's method for output-feedback controller presented in \cite{ilka2022novel} to solve the stabilization problem with $N$ control inputs. Firstly, we introduce the following matrix
\begin{equation}\label{Amatrix}
    A_i^{(k)} = A - \sum_{j=i+1}^{N} B_j R_{jj}^{-1} B_j^\top X_j^{(k)} - \sum_{j=1}^{i-1} B_j R_{jj}^{-1} B_j^\top X_j^{(k+1)},
\end{equation}
where $X_j^{(0)} = P_j$ for $j\in\mathcal{N}$ and $k$ is the iteration index. Then, the modified Newton's method, extended for the case with $N$ control inputs, is an asynchronous solving of LQR problem as in \cite{ilka2022novel}, but with varying plant matrix $A_i$ given in \eqref{Amatrix}. The algorithm for stabilization is given below. 
\begin{algorithm}
\caption{Modified Newton's Method For Stabilization of $N$-input system}
\label{MNMalg}
    \begin{enumerate}
    \item Initialize $R_{ii} > 0$ for $i\in\mathcal{N}$ and choose a small constant $\epsilon > 0$. Set the initial $X_i^{(0)} = P_i$ with $P_i$ computed via \eqref{lyapiter}. Set $k=0$ and $i=1$.
    \item Compute $A_i^{(k)}$ from \eqref{Amatrix}, matrix
    \begin{equation}\label{valuesMNM}
        G_i^{(k)} = R_{ii}^{-1}B_i^\top X_i^{(k)} (C_i^+ C_i - I_{n}),
    \end{equation}
    and 
    \begin{align}\label{step3}
         \begin{split}
                & \mathcal{R}(X_i^{(k)}) = C_i^\top C_i + G_i^{(k)\top} R_{ii} G_i^{(k)} + \\
                & A_i^{(k)\top} X_i^{(k)} + X_i^{(k)} A_i^{(k)} - X_i^{(k)} B_i R_{ii}^{-1} B_i^\top X_i^{(k)}.
         \end{split}
    \end{align}
    \item Solve 
    \begin{align}\label{step4}
        \begin{split}
        & (A_i^{(k)} - B_i R_{ii}^{-1} B_i^\top X_i^{(k)})^\top P_i^{(k)} + \\
        & P_i^{(k)} (A_i^{(k)} - B_i R_{ii}^{-1} B_i^\top X_i^{(k)}) = -  \mathcal{R}(X_i^{(k)}).
        \end{split}
    \end{align}
    with respect to $P_i^{(k)}$, and update $X_i^{(k+1)} = X_i^{(k)} + P_i^{(k)}$.
    \item Compute $\text{trace}(\mathcal{R}(X_i^{(k)})^\top \mathcal{R}(X_i^{(k)}))$. If $i=N$, go to step 5. Otherwise, set $i=i+1$ and perform steps 2-4.
    \item If for every $i\in\mathcal{N}$ $\text{trace}(\mathcal{R}(X_i^{(k)})^\top \mathcal{R}(X_i^{(k)})) < \epsilon$, then stop. Otherwise, set $i=1$, $k=k+1$, and repeat steps 2-4.
	\end{enumerate}
\end{algorithm}

The convergence proof of Algorithm \ref{MNMalg} is omitted due to space limitations and its similarity to the proof for the single-input case presented in \cite{ilka2022novel}. Note that, considering \eqref{valuesMNM} and \eqref{step3}, \eqref{step4} can be rewritten as
\begin{align}
&(A_i^{(k)} - B_i R_{ii}^{-1} B_i^\top X_i^{(k)})^\top (P_i^{(k)} + X_i^{(k)}) + \nonumber\\
&(P_i^{(k)} + X_i^{(k)}) (A_i^{(k)} - B_i R_{ii}^{-1} B_i^\top X_i^{(k)}) = \nonumber\\
&- C_i^\top C_i - G_i^{(k)\top} R_{ii} G_i^{(k)} - X_i^{(k)\top} B_i R_{ii} B_i^\top X_i^{(k)} \leq 0. \nonumber
\end{align}
Hence, if $A_i^{(k)} - B_i R_{ii}^{-1} B_i^\top X_i^{(k)}$ is stable, then $A_i^{(k)} - B_i R_{ii}^{-1} B_i^\top (P_i^{(k)} + X_i^{(k)})$ is also stable. For more details, we refer the reader to \cite{ilka2022novel} and \cite{engwerda2007algorithms}.

\begin{remark}
    The computation of $X_i^{(1)}$ via Newton's methods, i.e., the first iteration of the algorithm, might be problematic in some cases. This is a known drawback of Newton's methods referred in literature as ``disastrous first Newton step''. There exist step-size control techniques to overcome this problem described in \cite{guo2000newton} and \cite{benner1998exact}. Another problem might affect convergence of the algorithm happens if $A_i - B_i R_{ii}^{-1} B_i^\top X_i^*$  where $X_i^* = \lim_{k\to\infty} X_i^{(k)}$ has imaginary eigenvalues. This might affect significantly the convergence speed (for the details, we refer the reader to \cite{ilka2022novel,guo2000newton}). 
\end{remark}

After $(X_1,\dots,X_N)$ and associated $(K_1,\dots,K_N)$ as $K_i = R_{ii}^{-1} B_i^\top X_i C_i^+$, we apply the correction function \eqref{corrfun}, i.e.,  $\mathcal{G}(X_i, K_i,\alpha_i)$ to obtain $\tilde{X}_i$ such that   
\begin{equation}
    K_i C_i = R_{ii}^{-1} B_i^\top \tilde{X}_i = R_{ii}^{-1} B_i^\top \tilde{X}_i C_i^+ C_i.
\end{equation}
The obtained $\tilde{X}_i$ are the initialized stabilizing solutions, i.e., $X_i^{(0)} = \tilde{X}_i$ for $i\in\mathcal{N}$ are such that $A - \sum_{j=1}^N B_j R_{jj}^{-1} X_j^{(0)}$ is stable.

\subsection{Evaluation of the Parameters}
\label{sec:Evaluation of the Parameters}
In this section, we solve the inverse problem by finding $X_i$ such that $R_{ii}^{-1} B_i^\top X_i = K_{i,d} C_i$ and then substitute it in \eqref{ARE} to find $Q_i$ for $i\in\mathcal{N}$. To find such $X_i$, the gradient descent optimization is used again. Firstly, $R_{ij} = R_{ij}^\top$ such that $R_{ii} > 0$ for $i,j\in\mathcal{N}$ are initialized in the previous part. These parameters remain unchanged during the iterative procedure. To perform the optimization, we define
\begin{equation}
    e_i^{(s)} \coloneqq K_i^{(s)} C_i C_i^+ - K_{i,d} = R_{ii}^{-1} B_i^\top X_i^{(s)} C_i^+ - K_{i,d},
\end{equation}
which is a function of $X_i$ showing the difference between the desired $K_{i,d}$ and the feedback law calculated with $X_i^{(s)}$ where $s = 0,1,\dots$ is the iteration counter. Then, we construct the function $ e_i^{(s)\top} e_i^{(s)}$ which is convex in $X_i$ and has its minimum when $e_i^{(s)} = 0$, i.e., when $X_i^{(s)}$ is such that $K_i^{(s)} = R_{ii}^{-1} B_i^\top X_i^{(s)} C_i^+ = K_{i,d}$. Then, one can find the minimizing $X_i^{(s)}$ in the following way 
\begin{align}\label{gradupd}
    \begin{split}
        & e_i^{(s)} = R_{ii}^{-1} B_i^\top X_i^{(s)}  C_i^+ - K_{i,d},\\
        & X_i^{(s+1)}  = X_i^{(s)} - \beta_i \frac{\partial (\tr (e_i^{(s)\top} e_i^{(s)})}{\partial X_i},
    \end{split}
\end{align}
where $\beta_i > 0$ is the step size and the partial derivative is given by 
\begin{equation}\label{errtrace}
    \frac{\partial (\tr (e_i^{(s)\top} e_i^{(s)})}{\partial X_i} =  C_i^{+\top} B_i R_{ii}^{-1} e_i^{(s)} + e_i^{(s)\top} R_{ii}^{-1} B_i^\top C_i^+.
\end{equation}
At each iteration $s+1$, to ensure that there exists $K_i^{(s+1)}$ such that $K_i^{(s+1)} C_i =  R_{ii}^{-1} B_i^\top X_i^{(s+1)}$, we apply the $\mathcal{G}$ function, i.e., we compute corrected $\tilde{X}_i^{(s+1)} = \mathcal{G}(X_i^{(s+1)}, K_i^{(s+1)},\alpha_i)$. Then, $\tilde{X}_i^{(s+1)}$ is substituted into \eqref{ARE} rewritten as 
\begin{align}\label{IOCupdate}
    \begin{split}
        & Q_i^{(s+1)} = - \sum_{j=1}^N (K_j^{(s+1)} C_j)^\top R_{ij} (K_j^{(s+1)} C_j) - \\
        & (A - \sum_{j=1}^N B_j K_j^{(s+1)} C_j)^\top \tilde{X}^{(s+1)}_i - \\
        & \tilde{X}^{(s+1)}_i (A - \sum_{j=1}^N B_j K_j^{(s+1)} C_j). 
    \end{split}
\end{align}
to update $Q_i^{(s+1)}$ where $K_j^{(s+1)} C_j = R_{jj}^{-1} B_j^\top \tilde{X}_j^{(s+1)}$. Then $\tilde{X}_i^{(s+1)}$ is used as $X_i^{(s)}$ in \eqref{gradupd} for the next iteration. In that way, \eqref{gradupd} with the correction \eqref{corrfun} and \eqref{IOCupdate} are repeated until $\tilde{X}_i^{(s+1)}$ is such that $\|e_i^{(s+1)} e_i^{(s+1)\top}\| \leq \delta_i$ where $\delta_i > 0$ is a measure of the desired precision. Note that \eqref{gradupd} without \eqref{corrfun} already guarantees the convergence to the desired feedback law. Application of \eqref{corrfun} is only necessary to guarantee that $K_i^{(s+1)} C_i = R_{jj}^{-1} B_j^\top \tilde{X}_j^{(s+1)}$ holds at each iteration $s=0,1,\dots$. Also, using \eqref{IOCupdate}, one can estimate $Q_i^{(0)}$ using $X_i^{(0)}$ for $i\in\mathcal{N}$ if it is necessary.

Finally, we can summarize the whole procedure for the model-based algorithm that gives a stabilizing set of solutions is described below. 

\begin{algorithm}
\caption{Model-based Algorithm for Solving $N$-player Inverse Output-feedback Differential Game}
\label{mbalg}
    \begin{enumerate}
    \item Set $\alpha_i,\beta_i > 0 $ step sizes for $i,j \in\mathcal{N}$ and initialize $R_{ij} = R_{ij}^\top$ for $j\neq i$.
    \item Perform Algorithm \ref{MNMalg}, and derive the initial stabilizing tuple $(X_1,\dots, X_N)$.  
    \item Apply correction \eqref{corrfun} to  $(X_1,\dots, X_N)$ to receive $(X_1^{(0)},\dots, X_N^{(0)})$.
    \item Set $s=0$. Perform \eqref{gradupd} and then correction \eqref{corrfun} to $X_i^{(s+1)}$ to get $\tilde{X}_i^{(s+1)}$. 
    \item Update $Q_i^{(s+1)}$ as in \eqref{IOCupdate}.
    \item If $\|e_i^{(s+1)} e_i^{(s+1)\top}\| \leq \delta_i$ then stop. Otherwise, set $s=s+1$ and $X_i^{(s)} = \tilde{X}_i^{(s)}$, and repeat steps 4-5.
	\end{enumerate}
\end{algorithm}

Note that the correction procedure \eqref{corrfun} operates as an \textit{inner loop} within the outer loop responsible for driving $X_i^{(s)}$ to the desired value, where $i\in\mathcal{N}$. If the existence of the feedback law  (as well as the stability of the dynamics) at each iteration is not required, it is possible to use the gradient descent update to solve the inverse problem without applying the correction function \eqref{corrfun}. The way to do this is shown in Section \ref{sec:Distributed Procedure}.

\subsection{Analysis}
In this section, we introduce theoretical results associated to the model-based Algorithm \ref{mbalg}. We show its iterative stability and convergence to an NE solution. We also characterize the possible solutions for the inverse problem and show how to derive an infinite number of possible parameter combinations without the need to use the algorithm again.  

Note that $X_i$ in $K_i = R_{ii}^{-1} B_i^\top X_i C_i^+$ linearly affects $K_i$. Suppose we have $X_i^{(0)}$ in \eqref{gradupd} such that $K_i^{(0)} C_i = R_{ii}^{-1} B_i^\top X_i$ where $K_i^{(0)}$ exists for $i\in\mathcal{N}$, and $A - \sum_{j=1}^N B_j K_j^{(0)} C_j$ is stable. Performing both \eqref{gradupd} and then \eqref{corrfun}, one updates $(\tilde{X}_1^{(s+1)},\dots,\tilde{X}_N^{(s+1)})$ and can compute $(K_1^{(s+1)},\dots, K_N^{(s+1)})$ that change from stabilizing tuple $(K_1^{(0)},\dots,K_N^{(0)})$ to $(K_{1,d},\dots,K_{N,d})$ which is also stabilizing as it is stated in the problem formulation. Hence, the following can be proved. 

\begin{theorem}\label{stabth}
There exists a set $\{(\bar{\beta}_i\}_{i=1}^N$ of the step sizes in \eqref{gradupd} such that $(\tilde{X}_1^{(s+1)},\dots,\tilde{X}_N^{(s+1)})$ is stabilizable  for every $s=0,1,\dots$ iteration.
\end{theorem}

\textbf{Proof}. Suppose $(X_1^{(0)},\dots,X_N^{(0)})$ is stabilizable, i.e., $A - \sum_{i=1}^N B_i R_{ii}^{-1} B_i^\top X_i^{(0)}$ is stable.  Iterating \eqref{gradupd}, for each $s=0,1,\dots$, we get $X_i^{(s+1)}$ such that $K_i^{(s+1)}$ can be computed as $K_i^{(s+1)} =  R_{ii}^{-1} B_i^\top X_i^{(s+1)} C_i^+$, is driven closer to $K_{i,d}$. Because for each $i\in\mathcal{N}$ the feedback law $K_i^{(s)}$ is iterated from stabilizing $K_i^{(0)}$ and driven to stabilizing $K_{i,d}$, there exist $\{\bar{\beta}_i\}_{i=1}^N$ such that $(K_1^{(s+1)},\dots,K_N^{(s+1)})$ is stabilizable \cite{lancaster1995algebraic} ,\cite{Bertsekas/99}. Performing correction, we only change $X_i^{(s)}$ to $\tilde{X}_i^{(s)}$ such that $B_i^\top X_i^{(s+1)}  C_i^+ C_i = B_i^\top \tilde{X}_i^{(s+1)} = K_i^{(s+1)} C_i$. Thus, the stability property of $K_i^{(s+1)} = R_{ii}^{-1} B_i^\top \tilde{X}_i^{(s+1)} C_i^+$ does not change, i.e., $(K_1^{(s+1)},\dots, K_N^{(s+1)})$ is stabilizable. This completes the proof.

The following theoretical result on convergence can be deduced for the model-based algorithm. 
\begin{theorem}
    Given initialized $R_{ij} = R_{ij}^\top$ with $R_{ii} > 0$ for $i,j\in\mathcal{N}$, $\tilde{X}_i^{(s+1)}$ converges to $\tilde{X}_i^\infty$ and, as a result, $Q_i^{(s+1)}$ converges to $Q_i^\infty$ for $i\in\mathcal{N}$. For each iteration $s=1,2,\dots$, $(A,\mathbf{B},\mathbf{C},\mathbf{Q}^{(s)},\mathbf{R})$ as well as resulting $(A,\mathbf{B},\mathbf{C},\mathbf{Q}^\infty,\mathbf{R})$ has the output-feedback NE $(K_1^{(s)},\dots,K_N^{(s)})$ and $(K_{1,d},\dots,K_{N,d})$, respectively.
\end{theorem}

\textbf{Proof}. 
    Given $X_i^{(s)}$, $X_i^{(s+1)}$ received via \eqref{gradupd} is such that $e_i^{(s+1)} e_i^{(s+1)\top} < e_i^{(s)} e_i^{(s)\top}$ due to gradient descent optimization. Performing \eqref{corrfun}, one gets $\tilde{X}_i^{(s+1)}$ such that 
    \begin{equation}
        B_i^\top \tilde{X}_i^{(s+1)} = B_i^\top \tilde{X}_i^{(s+1)} C_i^+ C_i = B_i^\top X_i^{(s+1)} C_i^+ C_i. 
    \end{equation}
   Since $C_i$ is full row rank, $B_i^\top \tilde{X}_i^{(s+1)} C_i^+ = B_i^\top X_i^{(s+1)} C_i^+$ can be stated. Hence, $e_i^{(s+1)}$ does not change after the correction \eqref{corrfun}. Since $e_i^{(s)} e_i^{(s)\top}$ is a convex function with the minimum when $e_i^{(s)} = 0$, performing gradient descent update \eqref{gradupd}, we get
   $\lim_{s\to\infty} \tilde{X}_i^{(s)} = \tilde{X}_i^\infty$ such that $R_{ii}^{-1} B_i^\top \tilde{X}_i^\infty = K_{i,d} C_i$. Since $\tilde{X}_i^\infty$ converges, $\lim_{s\to\infty} Q_i^{(s)} = Q_i^\infty$ in \eqref{IOCupdate} can be deduced. This completes convergence part of the proof.

   Since $\tilde{X}_i^{(s)}$ is directly substituted into \eqref{IOCupdate}, $Q_i^{(s)}$ with the initialized $R_{ij}$ are such that $\tilde{X}_i^{(s)}$ is a solution of \eqref{ARE}. The correction \eqref{corrfun} is guaranteed to generate $\tilde{X}_i{(s)}$ such that $K_i^{(s)} C_i = R_{ii}^{-1} B_i^\top \tilde{X}_i^{(s)}$ where $K_i^{(s)}$ exists. Thus, considering Theorem \ref{stabth}, we can conclude that all the conditions of Theorem \ref{mainth} are satisfied. The same arguments work for $\tilde{X}_i^\infty, K_i^\infty, Q_i^\infty, R_{ij}$ for $i,j\in\mathcal{N}$. This completes the proof.

\begin{remark}
    Note that at each iteration $(\tilde{X}_1^{(s+1)},\dots, \tilde{X}_N^{(s+1)})$ is a unique solution of \eqref{IOCupdate} rewritten as a set of Lyapunov equations (with fixed $K_j^{(s+1)}$ for $j\in\mathcal{N}$)
    \begin{align}\label{lyapunoveq}
    \begin{split}
        & Q_i^{(s+1)} + \sum_{j=1}^N (K_j^{(s+1)} C_j)^\top R_{ij} K_j^{(s+1)} C_j + \\
        & (A - \sum_{j=1}^N B_j K_j^{(s+1)} C_j)^\top \tilde{X}^{(s+1)}_i + \\
        & \tilde{X}^{(s+1)}_i (A - \sum_{j=1}^N B_j K_j^{(s+1)} C_j) = 0. 
    \end{split}
\end{align}
because $A - \sum_{j=1}^N B_j K_j^{(s+1)} C_j$ is stable (if the condition of Theorem \ref{stabth} is satisfied) and 
\begin{equation}
     Q_i^{(s+1)} + \sum_{j=1}^N (K_j^{(s+1)} C_j)^\top R_{ij} K_j^{(s+1)} C_j
\end{equation}
is symmetric. 
\end{remark}

\begin{remark}
    Assumption \ref{CBassump} is necessary to guarantee that $C_i^{+\top} B_i = B_i^\top C_i^+ \neq \mathbf{0}_{p_i\times m_i}$. Otherwise, \eqref{errtrace} is equal to zero and, as a result, \eqref{gradupd} does not work. A less restrictive alternative to circumvent this issue is to introduce perturbations to either $B_i$ or $C_i$ during the gradient descent update \eqref{gradupd}. While this adjustment resolves the issue, it typically leads to slower convergence. Additionally, such an approach introduces estimation errors in the approximate solution.
    
    %For example, one can use $B_{i\nu}$ where $\nu$ is a small constant that perturbs those elements of $B_i$ that lead to $C_i B_i = 0$, i.e., for  %For example, one can use $B_{i\nu}$ where $\nu$ is a small constant that perturbs those elements of $B_i$ that lead to $C_i B_i = 0$, i.e., for 
    % \begin{equation}
    %     B_i = \begin{pmatrix}
    %         1 \\ 0
    %     \end{pmatrix}\quad\text{and}\quad C_i = \begin{pmatrix} 0 & 1\end{pmatrix},\quad\text{one uses}\quad   
    %     B_{i,\nu} = \begin{pmatrix}
    %         1 \\ \nu
    %     \end{pmatrix}.
    % \end{equation}

\end{remark}

Next, let us introduce the following notations
\begin{align}
    \begin{split}
       & A_{cl} = A - \sum_{j=1}^N B_j K_{j,d} C_j, \quad \Delta X_i = X_i - X_i^\prime,\\
       & \Delta Q_i = Q_i - Q_i^\prime,\quad R_{ij} = R_{ij} - R_{ij}^\prime.
    \end{split}
\end{align}
Suppose we found $X_i$ such that 
\begin{equation}
    K_i C_i = R_{ii}^{-1} B_i^\top X_i = R_{ii}^{-1} B_i^\top X_i C_i^+ C_i = K_{i,d} C_i,
\end{equation}
and computed $Q_i$ for some initialized $R_{ij}$ for $i,j\in\mathcal{N}$. Then, any $X_i^\prime$, $Q_i^\prime$ and $R_{ij}^\prime$ satisfying 
\begin{align}
        & A_{cl}^\top \Delta X_i + X_i A_{cl} + \Delta Q_i + \sum_{j=1}^N (K_{j,d} C_j)^\top \Delta R_{ij} (K_{j,d} C_j),\nonumber\\
        & (R_{ii}^\prime)^{-1} B_i^\top X_i^\prime = (R_{ii}^\prime)^{-1} B_i^\top X_i^\prime C_i^+ C_ i= K_{i,d} C_i,
\end{align}
is also a solution of \textbf{Problem 1}. The above can be verified by subtracting \eqref{ARE} associated with $X_i^\prime, Q_i^\prime, R_{ij}^\prime$ from \eqref{ARE} associated with $X_i, Q_i, R_{ij}$ for $i,j\in\mathcal{N}$. Setting $\Delta X_i = 0$ and $\Delta R_{ii} = 0$, one can deduce the following proposition.

\begin{proposition}
 Suppose the model-based algorithm generated $(A,\mathbf{B},\mathbf{C},\mathbf{Q}^\infty,\mathbf{R})$ with $(X_1^\infty,\dots,X_N^\infty)$ such that $K_i C_i = R_{ii}^{-1} B_i^\top X_i^\infty$ for $i\in\mathcal{N}$. Then, any $Q_i^\prime$ and $R_{ij}^\prime$ for $i,j\in\mathcal{N}$, $j\neq i$, satisfying
\begin{equation}
    (Q_i^\infty - Q_i^\prime) + \sum_{j=1,j\neq i}^N (K_{j,d} C_j)^\top (R_{ij} - R_{ij}^\prime) (K_{j,d} C_j) = 0,
\end{equation}
also form a game that has $(K_{1,d},\dots,K_{N,d})$ as an NE. 
\end{proposition}
 This comes from the fact that performing change of these parameters does not affect the solution of \eqref{lyapunoveq} or the feedback laws of players. Thus, the parameters can be recalculated without the need to use the algorithm again.

\section{Model-free Algorithm}\label{sec:ModelFree}
 In this section we develop a model-free algorithm. The model-free algorithm is an extension of the model-based one and, as a result, posses the same analytical properties. Later, we show how to modify the model-free algorithm to be used in a distributed manner. Although in a distributed way \textbf{Problem 2} can be solved, the distributed version of the algorithm does not posses all the properties of the model based. This is discussed in Section \ref{sec:Distributed Procedure}. 
 
The procedure is straightforward - we will modify each step introduced before for the model-based algorithm in a way to avoid using the system dynamics. 

The first step is the derivation of the initial stabilizing tuple $(X_1^{(0)},\dots,X_N^{(0)})$. Also notice that both the correction procedure \eqref{corrfun} and the gradient descent update \eqref{gradupd} involve the matrix $B_i$ for $i\in\mathcal{N}$. Thus, we need either to substitute it by something known or estimate it. As for the model-based algorithm, we start with initializing $R_{ij} = R_{ij}^\top$ with $R_{ii} > 0$ for $i,j\in\mathcal{N}$. To modify the equations in a way to avoid using the dynamics we exploit ideas presented in \cite{jiang2012computational} and \cite{vrabie_adaptive_2011} and referred to as integral reinforcement learning \cite{vamvoudakis2011online}.

Firstly, we modify \eqref{lyapiter} that allows to find a stabilizing tuple $(P_1,\dots,P_N)$. We write the original system as 
\begin{align}\label{disttraj}
    & \dot{x} = (A - \sum_{j=1}^N B_j K_{j,d} C_j) x + \sum_{j=1}^N B_j(K_{j,d} C_j x + u_j) = \nonumber\\
    & A x + \sum_{j=1}^N B_j u_j,
\end{align}
where $u_j = - K_{j,d} C_j x + \omega_j$ (since $F_j = K_{j,d} C_j$ as was set in the model-based algorithm), i.e., control inputs generated by the stabilizable feedback law $K_{j,d} C_j$ and disturbed by exploration noise $\omega_j$ (the need of noise is explained further). 

Then, we modify \eqref{lyapiter} by multiplying it by $x^\top$ and $x$ and changing the dynamics term in the following way
\begin{align}\label{modlyapiter}
        & (Ax - \sum_{j=i}^N B_j K_{j,d} C_j x -\sum_{j=1}^{i-1} B_j R_{jj}^{-1}B_j^\top P_j + \nonumber\\ 
        & \sum_{j=i}^N B_j(K_{j,d} C_j x + u_j) + \sum_{j=1}^{i-1} B_j (R_{jj}^{-1}B_j^\top P_j x + u_j))^\top P_i x + \nonumber\\
        & x^\top P_i (Ax - \sum_{j=i}^N B_j K_{j,d} C_j x -\sum_{j=1}^{i-1} B_j R_{jj}^{-1}B_j^\top P_j + \\ 
        & \sum_{j=i}^N B_j(K_{j,d} C_j x + u_j) + \sum_{j=1}^{i-1} B_j (R_{jj}^{-1}B_j^\top P_j x + u_j)) - \nonumber\\
        & 2 \Big(\sum_{j=i}^N B_j (K_{j,d} C_j x + u_j) + \sum_{j=1}^{i-1} B_j (R_{jj}^{-1}B_j^\top P_j + u_j)\Big)^\top P_i x = \nonumber\\
        & - x^\top (\tilde{Q}_i - (K_{i,d} C_i)^\top R_{ii} K_{i,d} C_i) x \nonumber
\end{align}

In \cite{jiang2012computational}, it is shown that if we integrate \eqref{modlyapiter} from $t$ to $t+\delta t$, then solving \eqref{lyapiter} is equivalent to solving
\begin{align}\label{mfinit}
    \begin{split}
        &x(t+\delta t)^\top P_i x(t+\delta t) - x(t)^\top P_i x(t) - \\
        & 2\int_{t}^{t+\delta t} \Big(\sum_{j=i}^N (u_j + K_{j,d} C_j x)^\top B_j^\top P_i x\,d\tau + \\
        & \sum_{j=1}^{i-1} (u_j + R_{jj}^{-1}B_j^\top P_j x)^\top B_j^\top P_i x\Big)\,d\tau=\\
        & -\int_{t}^{t+\delta t} x^\top (\tilde{Q}_i + (K_{i,d} C_i)^\top R_{ii} K_{i,d} C_i) x\,d\tau.
    \end{split}
\end{align}
Equation \eqref{mfinit} is solved with respect to $P_i$ and $B_j P_i$.  Note that in contrast to the model-based algorithm, instead of $C_i^\top C_i$, we set $\tilde{Q}_i = C_i^\top C_i + I_n > 0$ in the right side of the equation. It is done in order to guarantee that $P_i$ is positive definite (as shown in \cite{li1995lyapunov}) because  $C_i^\top C_i$ is not necessarily positive definite. Thus, we can compute $B_i$ for $i\in\mathcal{N}$ as $B_i = ((B_i^\top P_i) (P_i)^{-1})^\top$ after solving the system of equations \eqref{mfinit}. Hence, performing \eqref{mfinit} and computing some $P_i$ for $i\in\mathcal{N}$ such that $A - \sum_{j=1}^N B_j R_{jj}^{-1} B_j^\top P_j$ is stable, we derived stabilizable solutions necessary to apply Algorithm \ref{MNMalg}.

Next, we modify the steps of Algorithm \ref{MNMalg}. Equations \eqref{Amatrix} and \eqref{valuesMNM} in Step 2 are skipped because they are not required for the model-free implementation of Algorithm \ref{MNMalg}. Thus, since $B_i$ is estimated from the above for $i\in\mathcal{N}$, we only need to avoid using $A$ in \eqref{step3} and \eqref{step4}.

Firstly, we rewrite \eqref{step3} as
\begin{align}
\begin{split}
& \mathcal{R}(X_i^{(k)}) = C_i^\top C_i + (K_i^{(k)} C_i)^\top R_{ii} (K_i^{(k)} C_i) + \\
& (A_i^{(k)} - B_i K_i^{(k)} C_i)^\top X_i^{(k)} + X_i^{(k)} (A_i^{(k)} - B_i K_i^{(k)} C_i),
\end{split}
\end{align}
where $K_i^{(k)} = R_{ii}^{-1} B_i^\top X_i^{(k)} C_i^+$. Then, the above is modified in the same way as performed to get \eqref{mfinit} from \eqref{lyapiter}. This results in 
\begin{align}\label{mfstep3}
        & \mathcal{R}_x(X_i^{(k)}) = \nonumber\\
        & \int_t^{t+\delta t} x^\top(C_i^\top C_i + (K_i^{(k)} C_i)^\top R_{ii} (K_i^{(k)} C_i))x\,d\tau +  \nonumber\\
                &  x(t+\delta t)^\top X_i^{(k)} x(t+\delta t) - x(t)^\top X_i^{(k)} x(t) - \\
                & 2\int_t^{t+\delta t} \Big[\sum_{j=i+1}^{N}(u_j +  R_{jj}^{-1} B_j^\top X_j^{(k)} x)^\top B_j^\top +
                 \sum_{j=1}^{i-1} (u_j + \nonumber\\ 
                 & R_{jj}^{-1} B_j^\top X_j^{(k+1)} x)^\top B_j^\top + (u_i + K_i^{(k)} C_i x )^\top B_i^\top  \Big] X_i^{(k)} x \,d\tau,\nonumber
\end{align}
where for each $i\in\mathcal{N}$ $x$ is given as 
\begin{align}\label{traj}
    & \dot{x} = A x + \sum_{j=1}^N B_j u_j = B_i(K_i^{(k)} C_i x + u_i) + \\
    & (A - \sum_{j=i+1}^{N} B_j R_{jj}^{-1} B_j^\top X_j^{(k)} - \sum_{j=1}^{i-1} B_j R_{jj}^{-1} B_j^\top X_j^{(k+1)} - \nonumber\\
    & B_i K_i^{(k)} C_i) x + 
    \sum_{j=i+1}^{N} B_j (R_{jj}^{-1} B_j^\top X_j^{(k)} x +u_j) + \nonumber\\
    &\sum_{j=1}^{i-1} B_j (R_{jj}^{-1} B_j^\top X_j^{(k+1)} x + u_j) \nonumber
\end{align}
which is the same trajectory as \eqref{disttraj}. In the above, we again modified the dynamics and subtracted the changes afterwards. %Later, we show that the change of $\mathcal{R}(X_i^{(k)})$ to $\mathcal{R}_x (X_i^{(k)})$ is valid and that the noise terms arising from $u_j$ terms cancel out. 

Finally, we rewrite \eqref{step4} as was done with \eqref{modlyapiter} and \eqref{mfstep3}. Firstly, we modify the dynamics term in \eqref{step4} as
\begin{align}
        & x^\top(A  - \sum_{j=i}^{N} B_j R_{jj}^{-1} B_j^\top X_j^{(k)}  - \sum_{j=1}^{i-1} B_j R_{jj}^{-1} B_j^\top X_j^{(k+1)} )^\top \cdot \nonumber\\
        & \cdot P_i^{(k)} x + 
         x^\top P_i^{(k)} (A - \sum_{j=i}^{N} B_j R_{jj}^{-1} B_j^\top X_j^{(k)}  - \nonumber\\
        &\sum_{j=1}^{i-1}B_j R_{jj}^{-1} B_j^\top X_j^{(k+1)} )x 
         = -  x^\top \mathcal{R}_x (X_i^{(k)}) x.
\end{align}
Solving the above is equivalent to solving
    \begin{align}\label{mfstep4}
        & x(t+\delta t)^\top P_i^{(k)} x(t+\delta t) - x(t)^\top P_i^{(k)} x(t) - \nonumber\\
        & 2\int_t^{t+\delta t} \Big[\sum_{j=i}^{N}(u_j +  R_{jj}^{-1} B_j^\top X_j^{(k)} x)^\top B_j^\top +\nonumber\\
        & \sum_{j=1}^{i-1} (u_j + R_{jj}^{-1} B_j^\top X_j^{(k+1)} x)^\top B_j^\top\Big] P_i^{(k)} x \,d\tau= \nonumber\\
        & -  \mathcal{R}_x(X_i^{(k)}),
    \end{align}
along the following trajectory
\begin{align}
    & \dot{x} = (A - \sum_{j=i}^{N} B_j R_{jj}^{-1} B_j^\top X_j^{(k)} - \sum_{j=1}^{i-1} B_j R_{jj}^{-1} B_j^\top X_j^{(k+1)}) x +\nonumber \\ 
    & \sum_{j=i}^{N} B_j (R_{jj}^{-1} B_j^\top X_j^{(k)} x +u_j) +  \sum_{j=1}^{i-1} B_j (R_{jj}^{-1} B_j^\top X_j^{(k+1)} x + \nonumber\\ & u_j)  = A x + \sum_{j=1}^N B_j u_j,
\end{align}
which is the same as trajectories \eqref{disttraj} and \eqref{traj}.

Equation \eqref{mfstep4} is solved with respect to $P_i^{(k)}$. Note, \eqref{step3} and \eqref{step4} were both modified in the same way -- both multiplied by $x^\top$ and $x$, and then integrated from $t$ to $t+\delta t$.

After computing all the feedback laws $(K_1,\dots,K_N)$ and associated $(X_1,\dots,X_N)$ such that $K_i = R_{ii}^{-1} B_i^\top X_i C_i^+$, we apply the correction function \eqref{corrfun}, i.e.,  $\mathcal{G}(X_i, K_i,\alpha_i)$ to obtain $\tilde{X}_i$ such that   
\begin{equation}
    K_i C_i = R_{ii}^{-1} B_i^\top \tilde{X}_i = R_{ii}^{-1} B_i^\top \tilde{X}_i C_i^+ C_i,
\end{equation}
where $B_i$ are derived while solving \eqref{lyapiter}. The obtained $\tilde{X}_i$ are the initialized stabilizing solutions, i.e., $X_i^{(0)} = \tilde{X}_i$.

And finally, we change the stopping criteria for the model-free algorithm to be $\|R_x(X_i^{(k)})\| < \epsilon$ where $\epsilon$ is a small constant. Hence, model-free Algorithm \ref{MNMalg} replaces \eqref{step3} and \eqref{step4} by \eqref{mfstep3} and \eqref{mfstep4}, respectively. The model-free algorithm extension of Algorithm \ref{MNMalg} is given below.
\begin{algorithm}
\caption{Model-free Modified Newton's Method For Stabilization of $N$-input system}
\label{mfMNMalg}
    \begin{enumerate}
    \item Initialize $R_{ii} > 0$ for $i\in\mathcal{N}$ and choose a small constant $\epsilon > 0$. Set the initial $X_i^{(0)} = P_i$ with $P_i$ computed via \eqref{mfinit}. Set $k=0$ and $i=1$.
    \item Compute $\mathcal{R}_x(X_i^{(k)})$ as in \eqref{mfstep3}.
    \item Solve \eqref{mfstep4} with respect to $P_i^{(k)}$ and update $X_i^{(k+1)} = X_i^{(k)} + P_i^{(k)}$.
    \item Compute $\|R_x(X_i^{(k)})\| < \epsilon$. If $i=N$, go to step 5. Otherwise, set $i=i+1$ and perform steps 2-4.
    \item If for every $i\in\mathcal{N}$ $\|R_x(X_i^{(k)})\| < \epsilon$, then stop. Otherwise, set $i=1$, $k=k+1$, and repeat steps 2-4.
	\end{enumerate}
\end{algorithm}

Next, we modify the parameters evaluation part. \eqref{gradupd} does not require any modification. it involves $B_i$ which is computed in \eqref{mfinit} for every $i\in\mathcal{N}$. Hence, only \eqref{IOCupdate} requires modification which is done in the same way as in \eqref{modlyapiter}, i.e.,  
\begin{align}\label{IOCupdatemf}
    \begin{split}
        & \int_{t}^{t+\delta t} x^\top Q_i^{(s+1)} x\,d\tau = \\
        & - \int_{t}^{t+\delta t} x^\top \sum_{j=1}^N (K_j^{(s+1)} C_j)^\top R_{ij} (K_j^{(s+1)} C_j) x\,d\tau -\\
        & x^\top(t+\delta t) \tilde{X}_i^{(s+1)} x(t+\delta t) + x^\top(t) \tilde{X}_i^{(s+1)} x(t) + \\
        & 2 \int_{t}^{t+\delta t}  \sum_{j=1}^N  (u_j + K_j^{(s+1)} C_j x)^\top B_j^\top \tilde{X}_i^{(s+1)} x,
    \end{split}
\end{align}
where $K_j^{(s+1)} C_j = R_{jj}^{-1} B_j^\top \tilde{X}_j^{(s+1)}$.

The model-free version of Algorithm \ref{mbalg} is shown further as Algorithm \ref{mfalg}.
\begin{algorithm}
\caption{Model-free Algorithm for Solving $N$-player Inverse Output-feedback Differential Game}
\label{mfalg}
    \begin{enumerate}
    \item Set $\alpha_i,\beta_i > 0 $ step sizes for $i,j \in\mathcal{N}$ and initialize $R_{ij} = R_{ij}^\top$ for $j\neq i$.
    \item Perform Algorithm \ref{mfMNMalg}, and derive the initial stabilizing tuple $(X_1,\dots, X_N)$.
    \item Apply correction \eqref{corrfun} to  $(X_1,\dots, X_N)$ to receive $(X_1^{(0)},\dots, X_N^{(0)})$.
    \item Set $s=0$. Perform \eqref{gradupd} and then correction \eqref{corrfun} to $X_i^{(s+1)}$ to get $\tilde{X}_i^{(s+1)}$. 
    \item Update $Q_i^{(s+1)}$ as in \eqref{IOCupdatemf}.
    \item If $\|e_i^{(s+1)} e_i^{(s+1)\top}\| \leq \delta_i$ then stop. Otherwise, set $s=s+1$ and $X_i^{(s)} = \tilde{X}_i^{(s)}$, and repeat steps 4-5.
	\end{enumerate}
\end{algorithm}
\begin{remark}
    Overall, we need to generate a single trajectory 
    \begin{equation}
        A x + \sum_{j=1}^N B_j u_j,
    \end{equation}
    where $u_j = - K_{j,d} C_j x + \omega_j$. 
\end{remark}

\begin{remark}
The noise is injected because all the equations involving trajectories are solved using the batch least-square method \cite{devore2015probability}, and, hence, the persistence of excitation (PE) condition needs to be satisfied \cite{xue2021inverse}, \cite{modares__2015}, \cite{jiang2012computational}, \cite{vrabie_adaptive_2011}. The noise can be sinusoids of different frequencies or some random noise. We refer the reader to \cite{ioannou2006} for more details on that matter.  
\end{remark}
\begin{remark}
    Note that the noise does not induce any bias in the estimate of $P_i^{(k)}$ in \eqref{mfstep4} and $Q_i^{(s+1)}$ in \eqref{IOCupdatemf} for $i\in\mathcal{N}$. The equations are designed in such a way that if one performs the inverse procedure to every equation that involves the noisy trajectory using L'Hopital's rule \cite{rudin1976principles}, as performed in \cite{modares__2015}, one gets the same equations as in the model-based algorithm. Therefore, all the analytical properties presented for the model-based algorithm are valid for the model-free algorithm. This is shown below for \eqref{step4} and \eqref{mfstep4}
\end{remark}
To demonstrate the noise cancellation and equivalence of performing \eqref{step4} and \eqref{mfstep4}, we firstly work with \eqref{mfstep3} by dividing $R_x (X_i^{(k)})$ by $\delta t$ and take the limit to get 
\begin{align}
        & \lim_{\delta t \to 0} \mathcal{R}_x(X_i^{(k)}) = \\
                & \lim_{\delta t \to 0} \int_t^{t+\delta t} x^\top(C_i^\top C_i + (K_i^{(k)} C_i)^\top R_{ii} (K_i^{(k)} C_i))x\,d\tau + \nonumber\\
                & \lim_{\delta t \to 0} \Big( x(t+\delta t)^\top X_i^{(k)} x(t+\delta t) - x(t)^\top X_i^{(k)} x(t) \Big) - \nonumber\\
                & \lim_{\delta t\to 0} \Big(2\int_t^{t+\delta t} \Big[\sum_{j=i+1}^{N}(u_j +  R_{jj}^{-1} B_j^\top X_j^{(k)} x)^\top B_j^\top +
                 \sum_{j=1}^{i-1} (u_j + \nonumber\\ 
                 & R_{jj}^{-1} B_j^\top X_j^{(k+1)} x)^\top B_j^\top + (u_i + K_i^{(k)} C_i x )^\top B_i^\top  \Big] X_i^{(k)} x \,d\tau\Big)\nonumber.
\end{align}
By L'Hopital's rule \cite{rudin1976principles}, we get 
\begin{align}
        & \lim_{\delta t \to 0} \mathcal{R}_x(X_i^{(k)}) = x^\top(C_i^\top C_i + (K_i^{(k)} C_i)^\top R_{ii} (K_i^{(k)} C_i))x + \nonumber\\
        & \dot{x}^\top X_i^{(k)} x + x^\top X_i^{(k)} \dot{x} - \\
        & 2\Big[\sum_{j=i+1}^{N}(u_j +  R_{jj}^{-1} B_j^\top X_j^{(k)} x)^\top B_j^\top +
        \sum_{j=1}^{i-1} (u_j + \nonumber\\ 
        & R_{jj}^{-1} B_j^\top X_j^{(k+1)} x)^\top B_j^\top + (u_i + K_i^{(k)} C_i x )^\top B_i^\top  \Big] X_i^{(k)} x\nonumber.
\end{align}
Considering \eqref{disttraj} with $u_j = -K_{j,d} C_j x +\omega_j$ for $j\in\mathcal{N}$, we get 
\begin{align}
        & \lim_{\delta t \to 0} \mathcal{R}_x(X_i^{(k)}) = x^\top(C_i^\top C_i + (K_i^{(k)} C_i)^\top R_{ii} (K_i^{(k)} C_i))x + \nonumber\\
        & (A x + \sum_{j=1}^N B_j (-K_{j,d} C_j x +\omega_j))^\top X_i^{(k)} x + \nonumber\\
        & x^\top X_i^{(k)} (A x + \sum_{j=1}^N B_j (-K_{j,d} C_j x +\omega_j)) - \\
        & 2 \Big[\sum_{j=i+1}^{N}(-K_{j,d} C_j x +\omega_j + R_{jj}^{-1} B_j^\top X_j^{(k)} x)^\top B_j^\top +\nonumber\\
        & \sum_{j=1}^{i-1} (-K_{j,d} C_j x +\omega_j + R_{jj}^{-1} B_j^\top X_j^{(k+1)} x)^\top B_j^\top +\nonumber\\ 
        & (-K_{i,d} C_i x +\omega_i + K_i^{(k)} C_i x )^\top B_i^\top  \Big] X_i^{(k)} x\nonumber.
\end{align}
Then, terms with $\sum_{j=1} B_j (-K_{j,d} C_j x +\omega_j)$ that include the noise cancel out and we get 
\begin{align}
        & \lim_{\delta t \to 0} \mathcal{R}_x(X_i^{(k)}) = x^\top(C_i^\top C_i + (K_i^{(k)} C_i)^\top R_{ii} (K_i^{(k)} C_i))x + \nonumber\\
        & x^\top A^\top X_i^{(k)} x + x^\top X_i^{(k)} A x - 2 \Big[\sum_{j=i+1}^{N}(R_{jj}^{-1} B_j^\top X_j^{(k)} x)^\top B_j^\top \nonumber\\
        & + \sum_{j=1}^{i-1} (R_{jj}^{-1} B_j^\top X_j^{(k+1)} x)^\top B_j^\top + (K_i^{(k)} C_i x )^\top B_i^\top  \Big] X_i^{(k)} x,
\end{align}
which, considering \eqref{Amatrix} and \eqref{valuesMNM}, can be regrouped to receive the following equality
\begin{align}
    \begin{split}
        & \lim_{\delta t \to 0} \mathcal{R}_x(X_i^{(k)}) = x^\top ( C_i^\top C_i + G_i^{(k)\top} R_{ii} G_i^{(k)} + \\
        & A_i^{(k)\top} X_i^{(k)} + X_i^{(k)} A_i^{(k)} - X_i^{(k)} B_i R_{ii}^{-1} B_i^\top X_i^{(k)}) x = \\
        & x^\top \mathcal{R} (X_i^{(k)}) x.
    \end{split}
\end{align}
Further, we use the following notations
\begin{align}
        & \mathcal{D}(P_i^{(k)}) = (A_i^{(k)} - B_i R_{ii}^{-1} B_i^\top X_i^{(k)})^\top P_i^{(k)} + \nonumber\\
        & P_i^{(k)} (A_i^{(k)} - B_i R_{ii}^{-1} B_i^\top X_i^{(k)}),\\
        & \mathcal{D}_x (P_i^{(k)}) = x(t+\delta t)^\top P_i^{(k)} x(t+\delta t) - x(t)^\top P_i^{(k)} x(t) - \nonumber\\
        & 2\int_t^{t+\delta t} \Big[\sum_{j=i}^{N}(u_j +  R_{jj}^{-1} B_j^\top X_j^{(k)} x)^\top B_j^\top +\\
        & \sum_{j=1}^{i-1} (u_j +R_{jj}^{-1} B_j^\top X_j^{(k+1)} x)^\top B_j^\top\Big] P_i^{(k)} x \,d\tau\nonumber.
\end{align}
Next, we divide $\mathcal{D}_x(P_i^{(k)})$  by $\delta t$ and taking the limit to get 
    \begin{align}
        \begin{split}
        &  \lim_{\delta t \to 0} \mathcal{D}_x (P_i^{(k)}) = \\
        & \lim_{\delta t \to 0}  \Big(x(t+\delta t)^\top P_i^{(k)} x(t+\delta t) - x(t)^\top P_i^{(k)} x(t)\Big) - \\
        &  \lim_{\delta t \to 0}  2\int_t^{t+\delta t} \Big[\sum_{j=i}^{N}(u_j +  R_{jj}^{-1} B_j^\top X_j^{(k)} x)^\top B_j^\top +\\
        & \sum_{j=1}^{i-1} (u_j + R_{jj}^{-1} B_j^\top X_j^{(k+1)} x)^\top B_j^\top\Big] P_i^{(k)} x \,d\tau.
        \end{split}
    \end{align}
Using again L'Hopital's rule, we get 
    \begin{align}
        \begin{split}
        & \lim_{\delta t \to 0} \mathcal{D}_x (P_i^{(k)}) =  \Big( \dot{x}^\top P_i^{(k)} x + x^\top P_i^{(k)} \dot{x}\Big) - \\
        & 2\Big[\sum_{j=i}^{N}(u_j +  R_{jj}^{-1} B_j^\top X_j^{(k)} x)^\top B_j^\top +\\
        & \sum_{j=1}^{i-1} (u_j + R_{jj}^{-1} B_j^\top X_j^{(k+1)} x)^\top B_j^\top\Big] P_i^{(k)} x.
        \end{split}
    \end{align}
Considering \eqref{disttraj} with $u_j = -K_{j,d} C_j x +\omega_j$ for $j\in\mathcal{N}$, we can rearrange the terms to get 
    \begin{align}
        \begin{split}
        & \lim_{\delta t \to 0} \mathcal{D}_x (P_i^{(k)}) = x^\top \Big[(A_i^{(k)} - B_i R_{ii}^{-1} B_i^\top X_i^{(k)})^\top P_i^{(k)} + \\
        & P_i^{(k)} (A_i^{(k)} - B_i R_{ii}^{-1} B_i^\top X_i^{(k)})\Big] x = x^\top \mathcal{D} (P_i^{(k)}) x.
        \end{split}
    \end{align}
Thus, considering \eqref{step3} and \eqref{step4} in Algorithm \ref{MNMalg}, we write 
\begin{equation}
    \lim_{\delta t \to 0} ( \mathcal{D}_x (P_i^{(k)}) + \mathcal{R}_x(X_i^{(k)})) = x^\top (\mathcal{D} (P_i^{(k)}) + \mathcal{R}(X_i^{(k)}) x = 0.
\end{equation}
Hence, $P_i^{(k)}$ computed from \eqref{mfstep4} is such that \eqref{step4} is satisfied for $i\in\mathcal{N}$. 

The same procedure can be done to show solution equivalence of \eqref{lyapiter} and \eqref{mfinit}, and of \eqref{IOCupdate} and \eqref{IOCupdatemf}.
\begin{remark}
For the proposed algorithm, it is assumed that we can observe $x$. Essentially, it is desired that one can solve the problem solely by observing $y_i$. In \cite{li2011state}, it is demonstrated that given $y_{i,t}$, the observations of the outputs $y_i$, under certain conditions related to system observability, one can express this relationship as
\begin{equation}
x(t) = G_i y_{i,t}.
\end{equation}
Here, $G_i$ is a matrix dependent on $A, B_1, \dots, B_N, K_{1,d}, \dots, K_{N,d}, C_i$, and it cannot be directly estimated. In such cases, solving the inverse problem involves substituting $x(t)$ with $G_i^\prime y_{i,t}$, where $G_i^\prime = (G_i^\top G_i)^{-1} G_i^\top$. However, deriving the parameters in explicit form (i.e., $Q_i$ instead of $(G_i^\prime)^\top Q_i G_i^\prime$) in this scenario is not straightforward.

Most likely, either a completely different approach or significant modification of the presented algorithm is necessary if one aims to use only observations of the outputs to solve the inverse problem and obtain the parameters in explicit form. For further details, we refer the reader to \cite{li2011state}, \cite{modares2016optimal}, and \cite{xue2023inverse}.
\end{remark}

\subsection{Model-free Algorithm: Implementation}

In this section, we demonstrate how to perform computations in \eqref{mfinit}, \eqref{mfstep4}, and \eqref{IOCupdatemf}. This is accomplished by utilizing the following property of the Kronecker product
\begin{equation}\label{kronprop}
    (c^\top \otimes a^\top) \text{vec}(B) = a^\top B c,
\end{equation}
and the batch least-square method \cite{devore2015probability}.

Firstly, we work with \eqref{mfinit}. 
In accordance with \cite{jiang2012computational}, we introduce the following notations
\begin{align}
    \begin{split}
        \hat{D} &= [D_{11},2D_{12},\dots,2D_{1n}, D_{22},2D_{23},\dots,D_{nn}]^\top,\\
        \hat{x} &= [x_1^2,x_1 x_2, \dots, x_1 x_n, x_2^2, x_2 x_3,\dots, x_n^2]^\top,
    \end{split}
\end{align}
where $D_{l_1 l_2}$ is a particular element of a square matrix $D$ of dimension $n$, i.e., $(P_i)_{l_1l_2}$ for $l_1,l_2=1,\dots,n$. $\hat{D} \in\mathbb{R}^{n(n+1)/2}$ and $\hat{x}\in\mathbb{R}^{n(n+1)/2}$.

Using \eqref{kronprop}, we can rewrite terms of \eqref{mfinit} as 
\begin{align}
    \begin{split}
        & x^\top (t+\delta t) P_i x (t+\delta t) - x^\top (t) P_i x(t) =\\
        & (\hat{x}(t+\delta t) - \hat{x}(t)) \hat{P}_i, \\
        &  (u_j + K_{j,d} C_j x)^\top B_j^\top P_i x = ((x^\top\otimes u_j^\top) + \\
        & (x^\top\otimes x^\top)(I_n\otimes (K_{j,d} C_j)^\top))\text{vec}(B_j^\top P_i),\\
        & (u_j + R_{jj}^{-1}B_j^\top P_j x)^\top B_j^\top P_i x = ((x^\top\otimes u_j^\top) + \\
        & (x^\top\otimes x^\top)(I_n\otimes (R_{jj}^{-1}B_j^\top P_j)^\top))\text{vec}(B_j^\top P_i) \\
        & x^\top (\tilde{Q}_i + (K_{i,d} C_i)^\top R_{ii} K_{i,d} C_i) x = \\
        & (x^\top\otimes x^\top)\text{vec}(\tilde{Q}_i + (K_{i,d} C_i)^\top R_{ii} K_{i,d} C_i).\\
    \end{split}
\end{align}

In addition to the above, we define $\delta_{xx}$, $I_{xx}$ and $I_{xu_j}$
as 
\begin{align}\label{nottttations}
        &\delta_{xx}= [\hat{x}(t_1)-\hat{x}(t_0),\hat{x}(t_2) - \hat{x}(t_1),\dots, \hat{x}(t_s) - \hat{x}(t_{s-1})]^\top,\nonumber\\
        &I_{xx} = \\
        & \int_{t_0}^{t_1} (x\otimes x)\,d\tau,\int_{t_1}^{t_2} (x\otimes x)\,d\tau,\dots,\int_{t_{s-1}}^{t_s} (x\otimes x)\,d\tau]^\top, \nonumber\\
        &I_{xu_j} = \\
        & \int_{t_0}^{t_1} (x\otimes u_j)\,d\tau,\int_{t_1}^{t_2} (x\otimes u_j)\,d\tau,\dots,\int_{t_{s-1}}^{t_s} (x\otimes u_j)\,d\tau]^\top, \nonumber
\end{align}
where $0\leq t_{l-1}\leq t_{l}$ for $l\in\{0,1,\dots,s\}$. Although the data intervals do not need to be equal, in our simulation presented further, we use $t_l-t_{l-1} = \delta t$ for $l\in\{0,1,\dots,s\}$.

Then, \eqref{mfinit} can be rewritten as 
\begin{equation}\label{lq11}
    H_i\begin{pmatrix}
        \hat{P}_i\\ \text{vec}(B_1^\top P_i)\\ \dots \\ \text{vec}(B_i^\top P_i) \\ \dots \\\text{vec}(B_N^\top P_i)\\
    \end{pmatrix}= \Xi_i
\end{equation}
where 
\begin{align}
    \begin{split}
        H_i = &[\delta_{xx}, -2I_{xu_1} - I_{xx}(I_n\otimes (R_{11}^{-1}B_1^\top P_1)^\top), \dots, \\ 
        & -2I_{xu_i-1} - I_{xx}(I_n\otimes (R_{i-1,i-1}^{-1}B_{i-1}^\top P_{i-1})^\top), \\
        & -2I_{xu_i} - I_{xx}(I_n\otimes (K_{i,d} C_i)^\top),\dots,\\
        & -2I_{xu_N} - I_{xx}(I_n\otimes (K_{N,d} C_N)^\top)],\\
        \Xi_i= & -I_{xx} (\tilde{Q}_i + (K_{i,d} C_i)^\top R_{ii} K_{i,d} C_i).
    \end{split}
\end{align}
Then, \eqref{lq11} can be solved as 
\begin{equation}\label{lq111}
        \begin{pmatrix}
        \hat{P}_i\\ \text{vec}(B_1^\top P_i)\\ \dots \\ \text{vec}(B_i^\top P_i) \\ \dots \\\text{vec}(B_N^\top P_i)\\
    \end{pmatrix}= (H_i^\top H_i)^{-1} H_i^\top \Xi_i.
\end{equation}
After this is computed, one can estimate $B_i = ((B_i^\top P_i) (P_i)^{-1})^\top$ either for each $i$ from its associated \eqref{lq11} or from \eqref{lq11} associated with any player $i$ for all the players as $B_j = ((B_j^\top P_i) (P_i)^{-1})^\top$. Note that the vector of unknowns has $n(n+1)/2 + \sum_{j=1}^N m_j n$ parameters. Thus, we need enough data points to satisfy $s \geq n(n+1)/2 + \sum_{j=1}^N m_j n$. 
\begin{remark}
    If $H_i$ is an invertible square matrix, the right side of \eqref{lq111} can be computed as $H_i^{-1} \Xi_i$.
\end{remark}

Next, we work with \eqref{mfstep3}. Using the introduced above notations, we can rewrite \eqref{mfstep3} as 
\begin{align}
    \begin{split}
        & \mathcal{R}_x(X_i^{(k)}) = \\
        &\int_t^{t+\delta t} (x^\top\otimes x^\top) \text{vec} (C_i^\top C_i + (K_i^{(k)} C_i)^\top R_{ii} (K_i^{(k)} C_i))\,d\tau + \\
        & \hat{x}(t+\delta t) \hat{X}_i^{(k)} - \hat{x}(t) \hat{X}_i^{(k)} - 2\int_t^{t+\delta t} \Big[\sum_{j=i+1}^{N}
        ((x^\top\otimes u_j^\top) + \\
        & (x^\top\otimes x^\top)(I_n\otimes (R_{jj}^{-1} B_j^\top X_j^{(k)})^\top))\text{vec}(B_j^\top X_i^{(k)}) \\
        & \sum_{j=1}^{i-1} ((x^\top\otimes u_j^\top) + \\
        & (x^\top\otimes x^\top)(I_n\otimes (R_{jj}^{-1} B_j^\top X_j^{(k+1)})^\top))\text{vec}(B_j^\top X_i^{(k)}) + \\
        & ((x^\top\otimes u_i^\top) + \\
        & (x^\top\otimes x^\top)(I_n\otimes (K_i^{(k)} C_i)^\top))\text{vec}(B_j^\top X_i^{(k)})\Big] \,d\tau.
    \end{split}
\end{align}
Then, the collection of $\mathcal{R}_x(X_i^{(k)}) $ computed on intervals $t$ and $t+\delta t$ is denoted by $RX_i^{(k)}$ 
\begin{align}
        & RX_i^{(k)} = \\
        & I_{xx} (C_i^\top C_i + (K_i^{(k)} C_i)^\top R_{ii} (K_i^{(k)} C_i)) + \delta_{xx} \hat{X}_i^{(k)} - \nonumber\\ 
        & 2\Big[\sum_{j=i+1}^{N}
        (I_{xu_j} + I_{xx} (I_n\otimes (R_{jj}^{-1} B_j^\top X_j^{(k)})^\top))\text{vec}(B_j^\top X_i^{(k)}) \nonumber\\
        & \sum_{j=1}^{i-1} (I_{xu_j} + I_{xx}(I_n\otimes (R_{jj}^{-1} B_j^\top X_j^{(k+1)})^\top))\text{vec}(B_j^\top X_i^{(k)}) + \nonumber\\
        & (I_{xu_i} + I_{xx}(I_n\otimes (K_i^{(k)} C_i)^\top))\text{vec}(B_j^\top X_i^{(k)})\Big]\nonumber.
\end{align}
Next, equation \eqref{mfstep4} is rewritten as 
    \begin{align}\label{rxik}
        &  - RX_i^{(k)} = \delta_{xx} \hat{P}_i^{(k)} - \\ 
        & 2\Big[\sum_{j=i}^{N}
        (I_{xu_j} + I_{xx} (I_n\otimes (R_{jj}^{-1} B_j^\top X_j^{(k)})^\top))\text{vec}(B_j^\top P_i^{(k)}) +\nonumber\\
        & \sum_{j=1}^{i-1} (I_{xu_j} + I_{xx}(I_n\otimes (R_{jj}^{-1} B_j^\top X_j^{(k+1)})^\top))\text{vec}(B_j^\top P_i^{(k)})\Big]\nonumber.
    \end{align}
Then, using the notations 
\begin{align}
    \begin{split}
        & \Theta_i^{(k)} = [\delta_{xx}, I_{xu_1} + I_{xx}(I_n\otimes (R_{11}^{-1} B_1^\top X_1^{(k+1)})^\top), \dots, \\
        & I_{xu_{i-1}} + I_{xx}(I_n\otimes (R_{i-1,i-1}^{-1} B_{i-1}^\top X_{i-1}^{(k+1)})^\top), \\
        &I_{xu_{i}} + I_{xx} (I_n\otimes (R_{i,i}^{-1} B_{i}^\top X_{i}^{(k)})^\top),\dots, \\
        & I_{xu_N} + I_{xx} (I_n\otimes (R_{NN}^{-1} B_N^\top X_N^{(k)})^\top)]
    \end{split}
\end{align}
we can rewrite \eqref{rxik} as in \eqref{lq11} as follows
\begin{equation}\label{lq22}
    \Theta_i^{(k)} \begin{pmatrix}
        \hat{P}_i^{(k)} \\ \text{vec}(B_1^\top P_i^{(k)})\\ \dots \\ \text{vec}(B_i^\top P_i^{(k)}) \\ \dots \\\text{vec}(B_N^\top P_i^{(k)})\\
    \end{pmatrix}= RX_i^{(k)}.
\end{equation}
This can be solved as 
\begin{equation}\label{lq222}
        \begin{pmatrix}
        \hat{P}_i^{(k)} \\ \text{vec}(B_1^\top P_i^{(k)})\\ \dots \\ \text{vec}(B_i^\top P_i^{(k)}) \\ \dots \\\text{vec}(B_N^\top P_i^{(k)})\\
    \end{pmatrix}= (\Theta_i^{\top(k)} \Theta_i^{(k)})^{-1} \Theta_i^{(k)\top} RX_i^{(k)}.
\end{equation}
\begin{remark}
    If $\Theta_i^{(k)}$ is an invertible square matrix, the right side of \eqref{lq222} can be computed as $(\Theta_i^{(k)})^{-1} RX_i^{(k)}$.
\end{remark}
The dimension of the left-hand side is the same as in \eqref{lq111}. Thus, the same $s \geq n(n+1)/2 + \sum_{j=1}^N m_j n$ data intervals are required. 

The stopping criteria $\|R_x(X_i^{(k)})\| < \epsilon_i$ in Algorithm \ref{mfMNMalg} can be rewritten as 
\begin{equation}
    \|(I_{xx}^\top I_{xx})^{-1} I_{xx} RX_i^{(k)}\| < \epsilon_i.
\end{equation}

Finally, we show the implementation of \eqref{IOCupdatemf}. Using the same notations as in the above, we can rewrite \eqref{IOCupdatemf} as
\begin{align}
    \begin{split}
        & I_{qx} \hat{Q}_i^{(s+1)} = \\
        & - I_{xx} \sum_{j=1}^N (K_j^{(s+1)} C_j)^\top R_{ij} ((K_j^{(s+1)} C_j) -\delta_{xx} \hat{\tilde{X}}_i^{(s+1)} + \\
        & 2\sum_{j=1}^N (I_{xu_j}\text{vec}(B_j^\top \tilde{X}_j^{(s+1)}) + I_{xx} \text{vec}(K_j^{(s+1)} B_j^\top \tilde{X}_i^{(s+1)})) = \\
        & \Omega_i^{(s+1)},
    \end{split}
\end{align}
where $I_{qx}$ is 
\begin{equation}
           I_{qx} = [\int_{t_0}^{t_1} \hat{x}\,d\tau,\int_{t_1}^{t_2} \hat{x}\,d\tau,\dots,\int_{t_{s-1}}^{t_s} \hat{x}\,d\tau]^\top = \Omega_i^{(s+1)},
\end{equation}
Then, $\hat{Q}_i^{(s+1)}$ for any $s=0,1,\dots$ can be estimated as 
\begin{equation}\label{lq333}
    \hat{Q}_i^{(s+1)}= (I_{qx}^\top I_{qx})^{-1} I_{qx}^\top \Omega_i^{(s+1)}.
\end{equation}
The dimension of the left-hand side is the same as in \eqref{lq111}. Thus, the same $s \geq n(n+1)/2 + \sum_{j=1}^N m_j n$ data intervals are required. 

In fact, if one does not require $Q_i^{(s+1)}$ for every iteration $s$, the above does not need to be implemented interactively for every new $\tilde{X}_i^{(s+1)}$. Instead, one can implement it only once when the desired precision is reached, i.e., $\|e_i^{(s+1)} e_i^{(s+1)\top}\| \leq \delta_i$. This is done in the simulation example dedicated to the model-free algorithm (Section \ref{sec:MFSIM}).

\subsection{Model-free Algorithm: Distributed Procedure}
\label{sec:Distributed Procedure}

In this section, we suggest how the inverse problem can be solved in a distributed manner. Let us assume that there are $N$ agents, where each agent $a_i$ computes the cost function parameters of player $i\in\mathcal{N}$ and only knows $C_i$ and the target output-feedback law $K_{i,d}$; it can collect $x$ and $u_i$ by applying a control law to the system. Thus, agent $a_i$ knows the output matrix and can collect data associated only with player $i$. We provide the simplest scheme that solves the inverse problem, but the stability of the dynamics is not preserved during the iterative process.

We suggest that each agent $a_i$ initializes some $X_i^{(0)} = \mathbf{0}_n$. Then, each agent $a_i$ initializes some $R_{ii} > 0$, and at each iteration $s$ evaluates 
\begin{equation}
    e_i^{(s)} = R_{ii}^{-1} B_i^\top X_i^{(s)} - F_{i,d}, 
\end{equation}
where $F_{i,d} = K_{i,d} C_i$ for $i\in\mathcal{N}$, and performs the gradient descent update 
\begin{equation}\label{gradupddistr}
    X_i^{(s+1)} = X_i^{(s)} - \gamma_i \frac{\partial (e_i^{(s)\top}e_i^{(s)})}{\partial X_i}, 
\end{equation}
where $\gamma_i > 0$ is the step size and the partial derivative is given by 
\begin{equation}
    \frac{\partial (e_i^{(s)\top}e_i^{(s)})}{\partial X_i} = e_i^{(s)\top} R_{ii}^{-1} B_i^\top + B_i R_{ii}^{-1} e_i^{(s)}.
\end{equation}
Next, we need to estimate or substitute with something $B_i$. Consider
\begin{align}\label{distrlyap}
    \begin{split}
        & (A - \sum_{j=1}^N B_j F_{j,d})^\top P_i + P_i (A - \sum_{j=1}^N B_j F_{j,d}) =\\
        & - (\tilde{Q}_i- F_{i,d}^\top R_{ii} F_{i,d})
    \end{split}
\end{align}
and a specific trajectory
\begin{equation}\label{trajdistrib}
    \dot{x}_i = (A - \sum_{j=1}^N B_j K_{j,d} C_j) x_i + B_i (K_{i,d} C_i x + u_i) 
\end{equation}
where $u_i = -K_{i,d} C_i x_i + \omega_i$ where $\omega_i$ is the noise. Note that in this trajectory, only the control input of player $i$ is disturbed with noise in contrast to the trajectory presented in \eqref{disttraj}. Then \eqref{distrlyap} can be solved in the model-free fashion as 
\begin{align}\label{distr1}
    \begin{split}
        &x_i(t+\delta t)^\top P_i x_i(t+\delta t) - x_i(t)^\top P_i x_i(t) - \\
        & 2\int_{t}^{t+\delta t} (u_i + K_{i,d} C_i x_i)^\top B_i^\top P_i x_i\,d\tau =\\
        & -\int_{t}^{t+\delta t} x_i^\top (\tilde{Q}_i + (K_{i,d} C_i)^\top R_{ii} K_{i,d} C_i) x_i\,d\tau.
    \end{split}
\end{align}
This equation can be solved in the same way as \eqref{mfinit} but with respect to $P_i$ and $B_i^\top P_i$ which allows to estimate $B_i$. Note that we do not require any knowledge of player $j\in\mathcal{N}/\{i\}$. Then, the gradient descent update \eqref{gradupddistr} can be performed to estimate such a $X_i^*$ that for the initialized $R_{ii}$, $a_i$ has 
\begin{equation}
    R_{ii}^{-1} B_i^\top X_i^* = K_{i,d} C_i. 
\end{equation}
Note, we do not need Assumption \ref{CBassump} since no multiplication of $C_i^\top B_i$ is made.

Then, agent $a_i$ can set $R_{ij} = \mathbf{0}_{m_j}$ for $j\neq i$. In that way, performing evaluation of $Q_i^*$ as in \eqref{IOCupdatemf}, agent $a_i$ does not need any knowledge of the other players, i.e., 
\begin{align}\label{distr2}
    \begin{split}
        & \int_{t}^{t+\delta t} x_i^\top Q_i^* x_i\,d\tau = \\
        & - \int_{t}^{t+\delta t} x_i^\top (K_{i,d} C_j)^\top R_{ii} ((K_{i,d} C_i) x_i\,d\tau -\\
        & x_i^\top(t+\delta t) X_i^* x_i(t+\delta t) + x_i^\top(t) X_i^* x_i (t)- \\
        & 2 \int_{t}^{t+\delta t} (u_i + K_{i,d} C_i x_i)^\top B_i^\top X_i^* x_i,
    \end{split}
\end{align}
which is solved as \eqref{IOCupdatemf}.

Now, the main question is the following: How can each agent $a_i$ obtain the specific trajectories $(x_i,u_i)$ described in \eqref{trajdistrib}? We propose achieving this by introducing one additional agent, denoted as $a_0$, tasked solely with sending signals to other agents, indicating when to apply the control input. We envision $a_0$ situated at the center of a star-type network, connected to every agent, while each agent $a_i$ is only connected to $a_0$. Moreover, we assume that $a_0$, in a clockwise manner, signals to agent $a_i$ when to apply one of two control inputs: either $u_i = K_{i,d} C_i$ (noise-free) or $u_i = K_{i,d} C_i + \omega_i$ (noisy). This process can be outlined as follows
\begin{enumerate}
    \item $a_0$ sends a signal to every agent $a_i$ for $i=1,\dots,N$ to apply the feedback law $K_{i,d}$, resulting in the application of noise-free $u_i$ control input.
    \item $a_0$ sends a signal to agent $a_1$ to apply noisy $u_1$. Agent $a_1$ then collects $(x_1,u_1)$.
    \item $\dots$
    \item $a_0$ sends a signal to agent $a_N$ to apply noisy $u_N$. Agent $a_N$ collects $(x_N,u_N)$. 
\end{enumerate}
In this manner, each agent gathers all necessary data to perform \eqref{gradupddistr}, \eqref{distr1}, and \eqref{distr2}.

To summarize, to perform the distributed procedure, each agent $a_i$ needs to set $R_{ij} = \mathbf{0}_{m_j}$ for $j\neq i$ and collect the specific data $(x_i,u_i)$."
\begin{remark}
It is not necessary to associate each agent $a_i$ with player $i$. Suppose the number of players $N$ is even. Then we can have $N/2$ agents indexed by $j=1,\dots, N/2$, where each agent $a_j$ computes the parameters for players $i=2j-1$ and $i = 2j$ (i.e., $a_1$ estimates for players $i=1,2$, $a_2$ estimates for players $i=3,4$, etc.). Then, agent $a_j$ needs to know $(C_{2j-1}, C_{2j})$, $(K_{2j-1,d}, K_{2j,d})$, and collect the trajectory data for both players $i = 2j-1$ and $i = 2j$. Also, $R_{2j-1,2j}$ and $R_{2j,2j-1}$ can be set to non-zero matrices. The grouping can be done in any other way.
\end{remark}
\begin{remark}
In fact, it might be possible to perform Algorithm \ref{mfalg} in a distributed manner while keeping the system stable, as shown in Theorem \ref{stabth}. Only \eqref{IOCupdatemf} requires trajectories and knowledge of all the players, and this can be avoided by setting particular $R_{ij}$ matrices to zero. We believe that in such a case, one should allow agent $a_0$ to not only send signals to apply the control inputs but also to regulate the $\beta_i$ step sizes in \eqref{gradupd} and probably keep track of the execution of iterative procedures in Algorithm \ref{mfalg} to adjust $\beta_i$ used by agents $a_i$ for keeping the system stable.
\end{remark}

\section{Simulations}
\label{sec:Siumation}

In this section, we demonstrate the results of simulations for the model-based algorithm, model-free algorithm and the distributed procedure of the model-free algorithm. 

\subsection{Model-based Algorithm Simulation}
Consider the following continuous time system dynamics
\begin{equation}
	\dot{x} = A x + \sum_{i=1}^2 B_i u_i,\quad y_1 = C_1 x,\quad y_2 = C_2 x,
\end{equation}
where 
\begin{equation}
	A = \begin{pmatrix}
		1 & 1\\ 0 & 2 	\end{pmatrix},\quad B_1 = \begin{pmatrix}
		1\\ 0
	\end{pmatrix}, \quad B_2 = \begin{pmatrix}0 \\ 1\end{pmatrix},
\end{equation}
with $C_1 = B_1^\top$ and $C_2 = B_2^\top$.

\begin{figure} 
    \centering
       \includegraphics[width=1\linewidth]{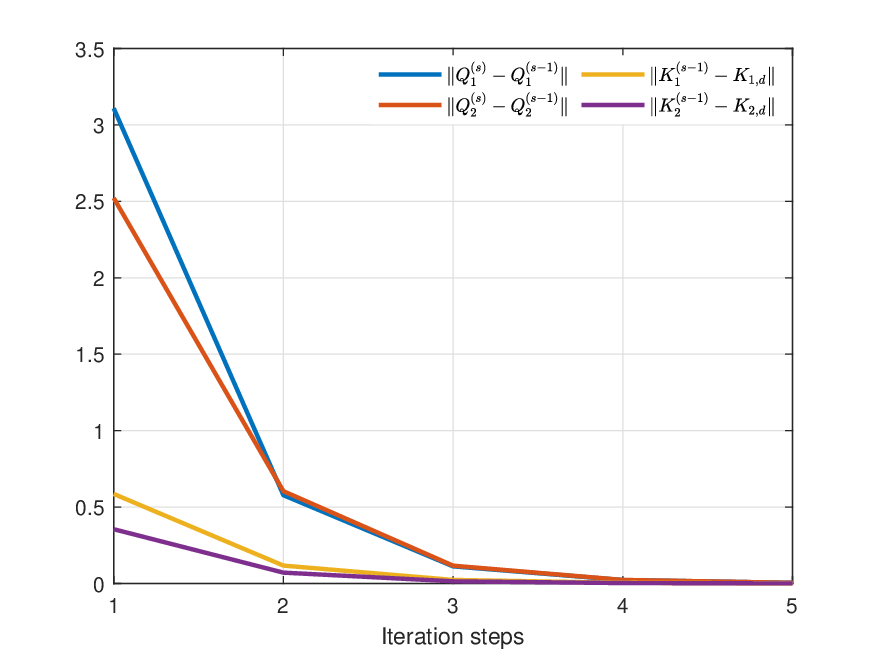}
  \caption{Algorithm \ref{mbalg}: convergence of the norm for iterations of $Q_i^{(s)}$ and $K_i^{(s)}$.}
  \label{fig1} 
\end{figure}

The target output-feedback laws are given by 
\begin{equation}
    K_{1,d} = 3,\quad K_{2,d} = 4. 
\end{equation}
As mentioned in Section \ref{sec:Initialization}, we initialized $R_{11} = R_{22} = 1$, 
% \begin{equation}
%     R_{11} = R_{22} = 1,
% \end{equation}
and computed $P_1$ and $P_2$ from \eqref{lyapiter} which gave us 
and the symmetric solution of AREs given by
\begin{align}
\begin{split}
    P_1 &= \begin{pmatrix}
          2.5  &  0.625 \\
    0.625  & 2.3125
    \end{pmatrix},\\ 
    P_2 &= \begin{pmatrix}
        1.0417  &  0.3348 \\
        0.3348  & 4.3616 
    \end{pmatrix}.
\end{split}
\end{align}
Then, Algorithm \ref{MNMalg} was applied to evaluate $K_1^{(0)} = 2.4142$ and $K_2^{(0)} = 4.3551$. After applying the correction \eqref{corrfun} (with the step sizes $\alpha_1 = 0.45$ and $\alpha_2 = 2 \alpha_1$), we got initial $X_1^{(0)}$ and $X_2^{(0)}$ as 
\begin{equation}
    X_1^{(0)} = \text{diag}(2.4142,8.3255),\, X_2^{(0)} = \text{diag}(1.0303,4.3551).
\end{equation}
Then, we set the step sizes for \eqref{gradupd} as $\beta_1 = \beta_2 = 0.6$ and $R_{12} = 2$, $R_{21} = 0.5$, and performed the procedure of Algorithm \ref{mbalg} to receive 
\begin{align}
\begin{split}
    Q_1^{(5)} &= \begin{pmatrix}
         3.0007 & -3 \\
   -3  & 3.6456 \\
    \end{pmatrix},\\ 
    Q_2^{(5)} &= \begin{pmatrix}
        -1.8 & -0.6751 \\
   -0.6751 & -0.0005
    \end{pmatrix}.
\end{split}
\end{align}
with $K_1^{(5)} = 3.0002$ and $K_2^{(5)} = 3.9999$,
and the symmetric stabilizing solution for \eqref{ARE} is
\begin{equation}
  X_1^{(5)} = \text{diag}(3.0002,8.9115),\, X_2^{(5)} = \text{diag}(0.6751,3.9999).  
\end{equation}
As the reader can see, the convergence is achieved fast -- $s=5$. The correction procedure \eqref{corrfun} is performed for $k=1,\dots,5$ at each iteration $s$. Figure \ref{fig1} demonstrates the convergence of the proposed algorithm. 
\begin{remark}
    The reader might notice that the step sizes $\alpha_1$ and $\alpha_2$ for players $1$ and $2$, respectively, are different. The reason is that the overshooting of the gradient descent method was observed. In fact, an adaptive learning rate might be used, e.g. Polyak step-size and the line search method \cite{sun2006optimization}.
\end{remark}

\subsection{Model-free Algorithm Simulation}
\label{sec:MFSIM}
\begin{figure} 
    \centering
  \subfloat[\label{fig2}]{%
        \includegraphics[width=0.5\linewidth]{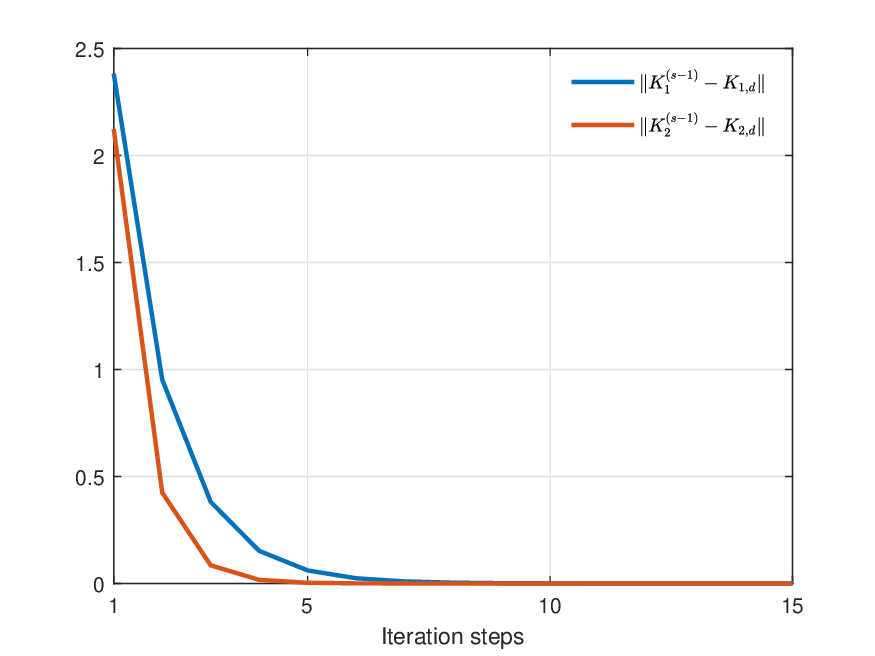}}
    \hfill
  \subfloat[\label{fig3}]{%
        \includegraphics[width=0.5\linewidth]{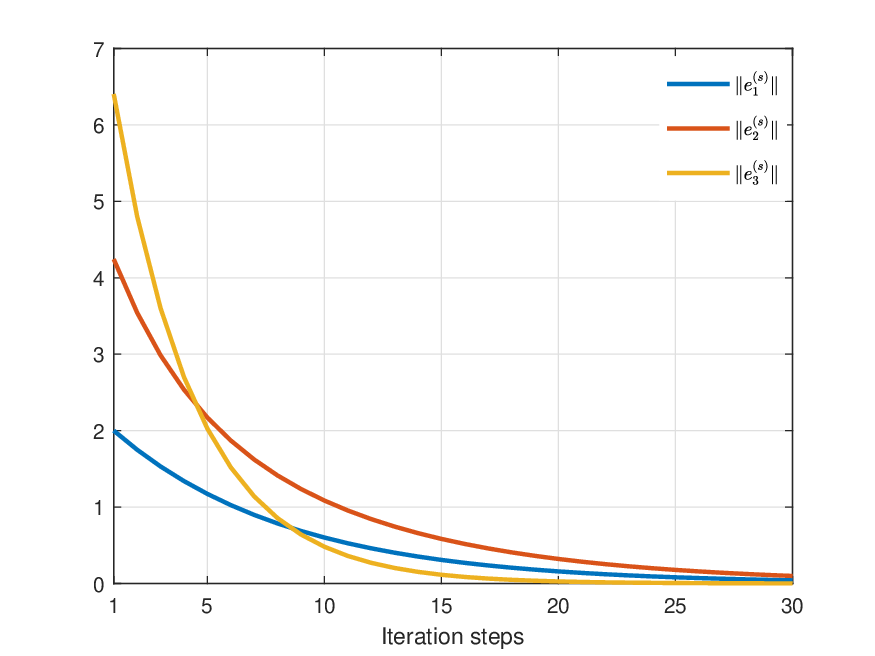}}
  \caption{(a) Algorithm \ref{mfalg}: convergence of the norm for iterations of $K_i^{(s)}$; Distributed procedure described in Section \ref{sec:Distributed Procedure}: convergence of the norm of $e_i^{(s)}$ in \eqref{gradupddistr}.}
\end{figure}

Consider the following continuous time system dynamics
\begin{align}
    \begin{split}
	& \dot{x} = A x + \sum_{i=1}^2 B_i u_i, \\
        & y_1 = C_1 x,\quad y_2 = C_2 x,
    \end{split}
\end{align}
where 
\begin{equation}
    A = \text{diag}(-3,0.5,4),\,B_1 = \begin{pmatrix}
		0 & 1 & 0
	\end{pmatrix}^\top,\, B_2 = \begin{pmatrix}0 & 0 & 1\end{pmatrix}^\top
\end{equation}
with $C_1 = B_1^\top$ and $C_2 = B_2^\top$.\\
The target output-feedback laws are given by 
\begin{equation}
    K_{1,d} = 4,\quad K_{2,d} = 6. 
\end{equation}
Firstly, we initialized $R_{11} = R_{22} = 1$. The noise used in this simulation is the same as in \cite{jiang2012computational}. For every player $i\in\mathcal{N}$, we set 
\begin{equation}\label{noiseused}
    \omega_i = 100 \sum_{j=1}^{1000} \sin(r_j t),
\end{equation}
where $r_j$ for $j=0,\dots,100$ is a randomly selected scalar from the range $[-500,500]$.\\ 
Next, we computed $P_1$ and $P_2$ from \eqref{modlyapiter} which gave us the symmetric solution of AREs given by
\begin{align}
    \begin{split}
         & P_1 = \text{diag}(0.1667,2.5714, 0.25),\\
         & P_2 = \text{diag}(0.1667,0.2414, 9.5).
    \end{split}
\end{align}
Then, Algorithm \ref{mfMNMalg} was applied to evaluate $K_1^{(0)} = 1.618$ and $K_2^{(0)} = 8.1231$.
After applying the correction \eqref{corrfun} (with the step sizes $\alpha_1 = 0.45$ and $\alpha_2 = 2 \alpha_1$), we got initial $X_1^{(0)}$ and $X_2^{(0)}$ as 
\begin{equation}
    X_1^{(0)} = \text{diag}(0,1.618,0),\quad X_2^{(0)} = \text{diag}(0,0,8.1231).
\end{equation}
Then, we set the step sizes for \eqref{gradupd} as $\beta_1 = \beta_2 = 0.6$ and $R_{12} = R_{21} = 0$, and performed the procedure of Algorithm \ref{mbalg} to receive  $K_1^{(15)} = 4$ and $K_2^{(15)} = 6$,
and the symmetric stabilizing solution for \eqref{ARE} is
\begin{align}
\begin{split}
    X_1^{(15)} &= \begin{pmatrix}
    2.382 & 0 &    2.382 \\
  0 &  4  &  0 \\
     2.382 & 0 &    2.382
    \end{pmatrix},\\ 
    X_2^{(15)} &= \begin{pmatrix}   
  -2.1231  & -2.1231  & 0 \\
   -2.1231 &  -2.1231  &  0 \\
    0 & 0 &   6
    \end{pmatrix}.
\end{split}
\end{align}
We did not performed computation of $Q_i^{(s+1)}$ at each iteration for $i=1,2$. Instead, we performed it once when the desired output-feedback laws were achieved. The results are
\begin{align}
\begin{split}
    Q_1^{(15)} &= \begin{pmatrix}
    14.2949 &   -0.0014  &  11.9104 \\
    -0.0014 & 12.0006 & -0.0002 \\
    11.9104 &  -0.0002 &  9.5279
    \end{pmatrix},\\ 
    Q_2^{(15)} &= \begin{pmatrix}
  -12.7386 & -13.8002  & 0 \\
  -13.8002 & -14.8617   & 0 \\
  0 & 0 & -12\\
    \end{pmatrix}.
\end{split}
\end{align}
The number of iterations is $s=15$. The correction procedure \eqref{corrfun} is performed for $k=1,\dots,15$ at each iteration $s$. Figure \ref{fig2} demonstrates the convergence of the proposed algorithm. 
% \begin{figure} 
%     \centering
%        \includegraphics[width=0.84\linewidth]{pic2.eps}
%   \caption{Algorithm \ref{mfalg}: convergence of the norm for iterations of $K_i^{(s)}$.}
%   \label{fig2} 
% \end{figure}

\subsection{Model-based Algorithm Simulation: Distributed Procedure}
\label{sec:MFSIMDP}

Consider the following continuous time system dynamics
\begin{equation}
	\dot{x} = A x + \sum_{i=1}^3 B_i u_i, 
\end{equation}
where 
\begin{equation}
	A = \begin{pmatrix}
		1 & 0 & 4\\ 0 & -3 & 2 \\ 0 & 5 & -2	\end{pmatrix}
\end{equation}
with 
\begin{align}
    \begin{split}
        & B_1 = (1,0,0)^\top, \quad B_2 = (0,1,0)^\top, \quad B_3 = (0,0,1)^\top, \\
        & C_1 = B_1^\top, \quad C_2 = (0,1,1), \quad C_3 = \begin{pmatrix}
            1 & 0 & 0 \\
            0 & 1 & 0
        \end{pmatrix}
    \end{split}
\end{align}
The target output-feedback laws are  
\begin{equation}
    K_{1,d} = 2,\quad K_{2,d} = 3,\quad K_{3,d} = \begin{pmatrix}4&5\end{pmatrix}.
\end{equation}
The initialized parameters $R_{ii}$ for $i=1,2,3$ are
\begin{equation}
    R_{11} = 4,\quad R_{22} = 3,\quad R_{33} = 2.
\end{equation}
The step sizes in \eqref{gradupddistr} are set to $\gamma_i = 1$ for $i=1,2,3$. The noise is in the previous simulation, i.e., \eqref{noiseused}.

After evaluation of $B_i$ for $i=1,2,3$ via \eqref{distr1}, the gradient descent optimization is applied to receive
\begin{align}
    \begin{split}
        & K_1^{(30)} = 1.9636, \quad K_2^{(30)} = 2.9554, \\
        & K_3^{(30)} = (3.9993, 4.9991),
    \end{split}
\end{align}
and the symmetric stabilizing solution for \eqref{ARE} is
\begin{align}
\begin{split}
    X_1^* = X_1^{(30)} &= \text{diag}(7.8543,0,0),\\
    X_2^* = X_2^{(30)} &= \begin{pmatrix}
           0 & 0 & 0 \\
   0 &    8.9952  &   8.7372 \\
   0 &    8.7372  & 0
    \end{pmatrix},\\ 
    X_3^* = X_3^{(30)} &= \begin{pmatrix}
        0 & 0 &  7.9986 \\
    0 & 0 &    9.9982 \\
    7.9986 &  9.9982  &  0
    \end{pmatrix}.
\end{split}
\end{align}
The computation of $Q_i^*$ via \eqref{distr2} for $i=1,2,3$ gives
\begin{align}
\begin{split}
    Q_1^* &= \begin{pmatrix}
           -0.286   & 0 &  -31.4174 \\
   0 & 0 & 0 \\
  -31.4174 & 0 &    0.0001
    \end{pmatrix},\\ 
    Q_2^* &= \begin{pmatrix}
           -0 &   34.9487 & 0  \\
   34.9487  & 80.9426 &  51.893 \\
   0  & 51.893  & -9.5026
    \end{pmatrix},\\ 
    Q_3^* &= \begin{pmatrix}
           31.988  & -0.0074 &  23.9957 \\
   -0.0074 &  -49.9999  & 79.9857 \\
   23.9957 &  79.9857 & -43.9921
    \end{pmatrix}.
\end{split}
\end{align}
The number of iterations is $s=30$. Figure \ref{fig3} demonstrates the convergence of the proposed algorithm. 
% \begin{figure} 
%     \centering
%        \includegraphics[width=0.84\linewidth]{pic3.eps}
%   \caption{Distributed procedure described in Section \ref{sec:Distributed Procedure}: convergence of the norm of $e_i^{(s)}$ in \eqref{gradupddistr}.}
%   \label{fig3} 
% \end{figure}

\section{Conclusion}
\label{sec:Conclusion}
In this paper, algorithms solving the inverse problem for LQ continuous time dynamic non-cooperative games are established. Firstly, we introduced the model-based algorithm and described its analytical properties. Then, the model-based algorithm is further extended to the model-free version that can solve the inverse problem in the case when the dynamics of the systems are unknown. In the process of development of the model-free algorithm, we modified algorithm in \cite{ilka2022novel} to work in a model-free fashion. Also, we introduced the distributed procedure for solving the inverse differential game problem. It is shown that all three algorithms generate a set of cost function parameters that form an equivalent game. It is shown how the resulted game can be modified without need to use the algorithm again.

In our work we considered the closed-loop infinite horizon games with linear dynamics. Thus, finite horizon games, games with non-linear dynamics and open loop games are interesting to study in the context of the inverse differential games with output-feedback. As mentioned before, finding the way to solve the inverse problem observing the output data $y_i$ of each player $i\in\mathcal{N}$ is a way for future research. Another option is to consider games with additional to stabilizability restrictions on the control inputs generated via output-feedback information structure. Essentially, differential games with stochastic elements in the dynamics are also of interest. 

\bibliographystyle{IEEEtran}
\bibliography{lat_temp.bib}

\end{document}